# PHASE TRANSITIONS FOR THE LONG-TIME BEHAVIOR OF INTERACTING DIFFUSIONS[1]

By A. Greven and F. den Hollander

*Universität Erlangen–Nürnberg and Leiden University*

Let $(\{X_i(t)\}_{i\in\mathbb{Z}^d})_{t\geq 0}$ be the system of interacting diffusions on $[0,\infty)$ defined by the following collection of coupled stochastic differential equations:

$$dX_i(t) = \sum_{j\in\mathbb{Z}^d} a(i,j)[X_j(t) - X_i(t)]\,dt + \sqrt{bX_i(t)^2}\,dW_i(t), \qquad i\in\mathbb{Z}^d, t\geq 0.$$

Here, $a(\cdot,\cdot)$ is an irreducible random walk transition kernel on $\mathbb{Z}^d\times\mathbb{Z}^d$, $b\in(0,\infty)$ is a diffusion parameter, and $(\{W_i(t)\}_{i\in\mathbb{Z}^d})_{t\geq 0}$ is a collection of independent standard Brownian motions on $\mathbb{R}$. The initial condition is chosen such that $\{X_i(0)\}_{i\in\mathbb{Z}^d}$ is a shift-invariant and shift-ergodic random field on $[0,\infty)$ with mean $\Theta\in(0,\infty)$ (the evolution preserves the mean).

We show that the long-time behavior of this system is the result of a delicate interplay between $a(\cdot,\cdot)$ and $b$, in contrast to systems where the diffusion function is subquadratic. In particular, let $\widehat{a}(i,j)=\frac{1}{2}[a(i,j)+a(j,i)]$, $i,j\in\mathbb{Z}^d$, denote the symmetrized transition kernel. We show that:

(A) If $\widehat{a}(\cdot,\cdot)$ is recurrent, then for any $b>0$ the system locally dies out.

(B) If $\widehat{a}(\cdot,\cdot)$ is transient, then there exist $b_*\geq b_2>0$ such that:

(B1) The system converges to an equilibrium $\nu_\Theta$ (with mean $\Theta$) if $0<b<b_*$.

(B2) The system locally dies out if $b>b_*$.

(B3) $\nu_\Theta$ has a finite 2nd moment if and only if $0<b<b_2$.

---

Received January 2005; revised October 2006.

[1]Supported by DFG and NWO through the German Priority Program "Interacting Stochastic Systems of High Complexity" and the Dutch-German Bilateral Group "Mathematics of Random Spatial Models from Physics and Biology."

(B4) The 2nd moment diverges exponentially fast if and only if $b > b_2$.

The equilibrium $\nu_\Theta$ is shown to be associated and mixing for all $0 < b < b_*$. We argue in favor of the conjecture that $b_* > b_2$. We further conjecture that the system locally dies out at $b = b_*$.

For the case where $a(\cdot,\cdot)$ is symmetric and transient we further show that:

(C) There exists a sequence $b_2 \geq b_3 \geq b_4 \geq \cdots > 0$ such that:

(C1) $\nu_\Theta$ has a finite $m$th moment if and only if $0 < b < b_m$.

(C2) The $m$th moment diverges exponentially fast if and only if $b > b_m$.

(C3) $b_2 \leq (m-1)b_m < 2$.

(C4) $\lim_{m\to\infty}(m-1)b_m = c = \sup_{m\geq 2}(m-1)b_m$.

The proof of these results is based on self-duality and on a representation formula through which the moments of the components are related to exponential moments of the collision local time of random walks. Via large deviation theory, the latter lead to variational expressions for $b_*$ and the $b_m$'s, from which sharp bounds are deduced. The critical value $b_*$ arises from a stochastic representation of the Palm distribution of the system.

The special case where $a(\cdot,\cdot)$ is simple random walk is commonly referred to as the parabolic Anderson model with Brownian noise. This case was studied in the memoir by Carmona and Molchanov [*Parabolic Anderson Problem and Intermittency* (1994) Amer. Math. Soc., Providence, RI], where part of our results were already established.

## 1. Introduction and main results.

1.1. *Motivation and background.* This paper is concerned with the long-time behavior of a particular class of systems with interacting components. In this class, the components are interacting diffusions that take values in $[0,\infty)$ and that are labelled by a countably infinite Abelian group $I$. The reason for studying these systems is two-fold: *new phenomena* occur, and a number of methodological problems can be tackled that are unresolved in the broader context of interacting systems with noncompact components. We begin by describing in more detail the background of the questions to be addressed.

A large class of interacting systems has the property that single components change according to a certain random evolution, while the interaction between the components is linear and can be interpreted as migration of mass, charge or particles. Examples are:

(1) *Interacting particle systems*, for example, voter model [34], branching random walk [22, 36], generalized potlatch and smoothing process [35], binary path process [31], coupled branching process [28, 29], locally dependent branching process [3], catalytic branching [27, 37].

(2) *Interacting diffusions*, for example, Fisher–Wright diffusion [9, 10, 13, 14, 24, 25, 32, 33, 40, 41, 44], critical Ornstein–Uhlenbeck process [19, 20], Feller's branching diffusion [15, 41], parabolic Anderson model with Brownian noise [7].
(3) *Interacting measure-valued diffusions*, for example, Fleming–Viot process [17], mutually catalytic diffusions [18], catalytic interacting diffusions [30].

Most of these systems display the following *universality*: independently of the nature of the random evolution of single components, the ergodic behavior of the system depends only on recurrence versus transience of the migration mechanism. More precisely, if the symmetrized migration kernel is recurrent then the system approaches trivial equilibria (concentrated on the "traps" of the system), whereas if the symmetrized migration kernel is transient then nontrivial extremal equilibria exist that can be parametrized by the spatial density of the components.

In this paper we study an example in a *different universality class*, one where the nature of the random evolution of single components does influence in a crucial way the long-time behavior of the system. In particular, we consider a system where the components evolve as diffusions on $[0,\infty)$ with diffusion function $bx^2$ and interact linearly according to a random walk transition kernel. Such a system is called the *parabolic Anderson model with Brownian noise* in the special case where the random walk is simple. In the recurrent case the system, as before, approaches a trivial equilibrium (concentrated on the "trap" with all components 0), so local extinction prevails. However, in the transient case we find three regimes, separated by critical thresholds $b_* > b_2 > 0$ (see Figure 1):

(I) ("low noise") $0 < b < b_2$: equilibria with finite 2nd moment.
(II) ("moderate noise") $b_2 \le b < b_*$: equilibria with finite 1st moment and infinite 2nd moment.
(III) ("high noise") $b \ge b_*$: local extinction.

We will show that the strict inequality $b_* > b_2$ depends on a large deviation principle for a renewal process in a random environment. This large deviation principle will be addressed in a forthcoming paper (Birkner, Greven and Hollander [5]). Local extinction at $b = b_*$ is a subtle issue that remains open.

For the case where the random walk transition kernel is *symmetric* we do a finer analysis. We show that in regime (I) there exists a sequence $b_2 \ge b_3 \ge b_4 \ge \cdots > 0$ such that the equilibria have a finite $m$th moment if and only if $0 < b < b_m$, while the $m$th moment diverges exponentially fast if and only if $b > b_m$ (see Figure 1). Moreover, we show that $b_2 \le (m-1)b_m < 2$ and that $\lim_{m\to\infty}(m-1)b_m = c = \sup_{m\ge 2}(m-1)b_m$. We show that in regimes (I) and (II) the equilibria are associated and mixing. We show that the critical value $b_*$ separating regimes (II) and (III) is linked to the Palm distribution of the system.

The reason for the above phase diagram is that there are two *competing* mechanisms: the migration pushes the components toward the mean value of the initial configuration, while the diffusion pushes them toward the boundary of the state space. Hence, there is a dichotomy in that either the migration dominates (giving nontrivial equilibria) or the diffusion dominates (giving local extinction). In the class of interacting diffusions we are concerned with here, *the migration and the diffusion have a strength of the same order of magnitude* and therefore the precise value of the diffusion parameter in relation to the migration kernel is crucial for the ergodic behavior of the system.

Our results are a completion and a generalization of the results in the memoir of Carmona and Molchanov [7]. In [7], Chapter III, the focus is on the *annealed Lyapunov exponents* for simple random walk, that is, on $\chi_m(b)$, the exponential growth rate of the $m$th moment of $X_0(t)$, for successive $m$. It is shown that for each $m$ there is a critical value $b_m$ where $\chi_m(b)$ changes from being zero to being positive (see Figure 2), and that the sequence $(b_m)$ has the qualitative properties mentioned earlier, that is, $b_m = 0$ for all $m$ in $d = 1, 2$ (recurrent case) and $b_2 \geq b_3 \geq b_4 \geq \cdots > 0$ in $d \geq 3$ (transient case). No existence of and convergence to equilibria is established below $b_2$, nor is any information on the equilibria obtained. There is also no analysis of what happens at the critical values. In our paper we are able to handle these issues due to the fact that we have variational expressions for $\chi_m(b)$ and $b_m$, which give us better control. In addition, we are able to get sharp bounds on $b_m$ that are valid for arbitrary symmetric random walk, which results in strict inequalities between the first few $b_m$'s.

In [7], Chapter IV, an analysis is given of the *quenched Lyapunov exponent* for simple random walk, that is, on $\chi_*(b)$, the a.s. exponential growth rate of $X_0(t)$. It is shown that $\chi_*(b)$ is negative for all $b > 0$ in $d = 1, 2$ (recurrent case), negative for $b > b_*$ and zero for $0 < b \leq b_*$ in $d \geq 3$ (transient case) for some $b_* \geq b_2$ (see Figure 2). This corresponds to the crossover at $b_*$, except for the proof that $b_* > b_2$, which we defer to a forthcoming paper [5]. In [7], Chapter IV, it is further shown that $\chi_*(b)$ has a singular asymptotics for $b \to \infty$. This asymptotics has been sharpened in a sequence of subsequent papers by Carmona, Molchanov and Viens [8], Carmona, Koralov and Molchanov [6] and Cranston, Mountford and Shiga [12].

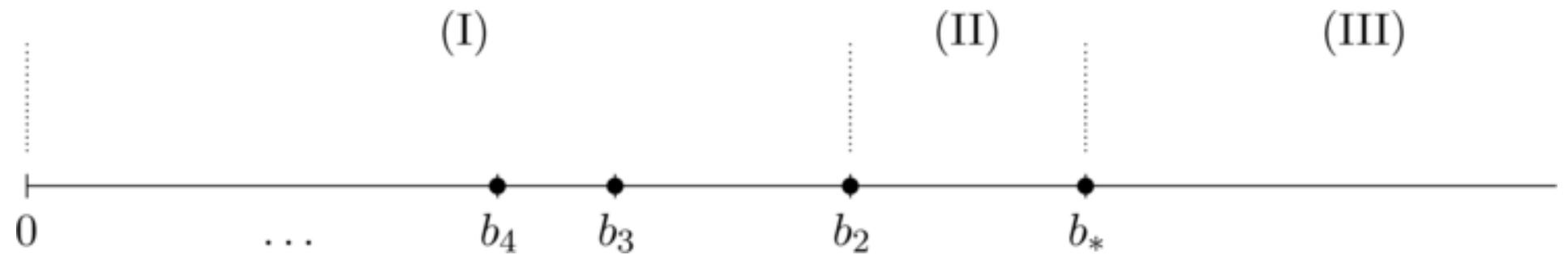

Fig. 1. *Phase diagram for the transient case.*

A scenario as described above is expected to hold for a number of interacting systems where the components take values in a noncompact state space, for example, generalized potlatch and smoothing [35] and coupled branching [28, 29]. But for none of these systems has the scenario actually been fully proven.

1.2. *Open problems.* We formulate a number of open problems that are not addressed in the present paper:

(A) Show that $b_2 > b_3 > b_4 > \cdots$. This property is claimed in [7], Chapter III, Section 1.6, but no proof is provided. We are able to show that $b_2 > b_3 > \cdots > b_m$ for an arbitrary symmetric random walk for which the average number of returns to the origin is $\leq 1/(m-2)$. For $m = 3$, this includes simple random walk in $d \geq 3$.

(B) Show that the system locally dies out at the critical value $b_*$.

(C) Show that $\chi_*(b) < 0$ for $b > b_*$, that is, show that there is no intermediate regime where the system locally dies out but only subexponentially fast. Shiga [41] has shown that the system locally dies out exponentially fast for $b$ sufficiently large.

(D) Find out whether there exists a characterization of $b_*$ in terms of the collision local time of random walks. This turns out to be a subtle problem, which has analogues in other models (see Birkner [3]). We find that such a characterization does exist for $b_m$ and for a certain $b_{**}$ with $b_* \geq b_{**}$. We have a characterization of $b_*$ in terms of the Palm distribution of our process, but this is relatively inaccessible. It therefore is a subtle problem to decide whether $b_* = b_{**}$ or $b_* > b_{**}$.

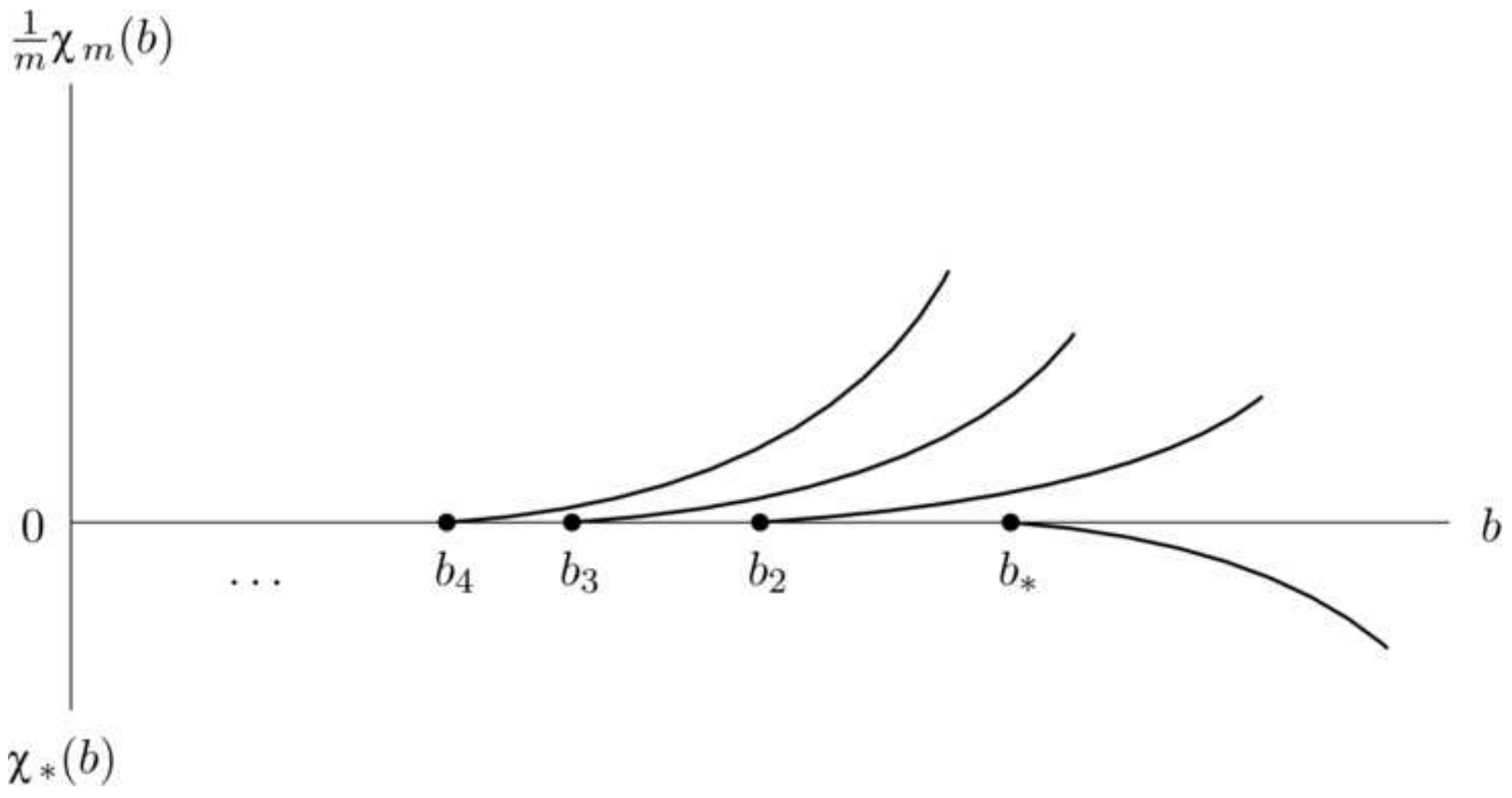

FIG. 2. *Qualitative picture of $b \mapsto \frac{1}{m}\chi_m(b)$ and $b \mapsto \chi_*(b)$ for the transient case.*

1.3. *Outline.* In Section 1.4 we define the model, formulate a theorem by Shiga and Shimizu [42] stating that our system of interacting diffusions has a unique strong solution with the Feller property, and introduce some key notions. In Section 1.5 we formulate two more theorems, due to Shiga [41] and to Cox and Greven [10], stating that our system locally dies out in the recurrent case and has associated mixing equilibria with finite 2nd moment in the transient case in regime (I). We complement these two theorems with two new results, stating that our system has associated mixing equilibria with finite 1st moment in the transient case in regime (II) and no equilibria in the transient case in regime (III). In Section 1.6 we present our finer results for regime (I), and have a closer look at regimes (II) and (III) as well, although much less detailed information is obtained for these regimes.

Sections 2–4 contain the proofs. Section 2 is devoted to moment calculations, which are based on a (Feynman–Kac type) representation formula for the solution of our system due to Shiga [41]. Through this representation formula, we express the moments of the components of our system in terms of exponential moments of the collision local time of random walks. Through the latter we are able to establish convergence to a (possibly trivial) equilibrium and to prove that this equilibrium is shift invariant, ergodic and associated. In Section 3 we study the exponential moments of the collision local time with the help of large deviation theory, which leads to a detailed analysis of the critical thresholds $b_m$ as a function of $m$ in regime (I), as well as to a description of the behavior of the system at $b_m$. Section 4 looks at survival versus extinction and relates the critical threshold $b_*$ between regimes (II) and (III) to the so-called Palm distribution of our system, where the law of the process is changed by size biasing with the value of the coordinate at the origin. There we argue in favor of the strict inequality $b_* > b_2$, which relies on an explicit representation formula for the Palm distribution.

1.4. *The model.* The models that we consider are systems of interacting diffusions $X = (X(t))_{t\geq 0}$, where

$$X(t) = \{X_i(t)\}_{i\in I} \in [0,\infty)^I, \tag{1.1}$$

with $I$ a countably infinite Abelian group. The evolution is defined by the following system of stochastic differential equations (SSDE):

$$dX_i(t) = \sum_{j\in I} a(i,j)[X_j(t) - X_i(t)]\,dt + \sqrt{bX_i(t)^2}\,dW_i(t), \qquad i\in I, t\geq 0. \tag{1.2}$$

Here:

(i) $a(\cdot,\cdot)$ is a Markov transition kernel on $I \times I$.
(ii) $b \in (0,\infty)$ is a parameter.
(iii) $W = (\{W_i(t)\}_{i\in I})_{t\geq 0}$ is a collection of independent standard Brownian motions on $\mathbb{R}$.

Equation (1.2) arises as the continuum limit of a self-catalyzing branching Markov chain whose branching rate depends on the local population size. As initial condition we take

$$X(0) \in \mathcal{E}_1, \tag{1.3}$$

where

$$\mathcal{E}_1 = \left\{ x = (x_i)_{i\in I} \in [0,\infty)^I : \sum_{i\in I} \gamma_i x_i < \infty \right\} \subset L^1(\gamma) \tag{1.4}$$

for any $\gamma = (\gamma_i)_{i\in I}$ satisfying the requirements

$$\begin{aligned} &\gamma_i > 0 \qquad \forall i \in I,\\ &\sum_{i\in I} \gamma_i < \infty,\\ &\exists M < \infty : \sum_{i\in I} \gamma_i a(i,j) \leq M\gamma_j \qquad \forall j \in I. \end{aligned} \tag{1.5}$$

We endow $\mathcal{E}_1$ with the product topology of $[0,\infty)^I$ (see Liggett and Spitzer [38]).

Since $|I| = \infty$, it is not possible to define the process uniquely in the strong sense on $[0,\infty)^I$ without putting growth conditions on the initial configuration, as in (1.4). However, the dependence of $\mathcal{E}_1$ on $\gamma$ is not very serious. For example, every probability measure $\rho$ on $[0,\infty)^I$ satisfying $\sup_{i\in I} E^\rho(X_i) < \infty$ is concentrated on $\mathcal{E}_1$ regardless of the $\gamma$ chosen ($E^\rho$ denotes expectation with respect to $\rho$). We also need the space $\mathcal{E}_2 \subset L^2(\gamma)$, which is defined as in (1.4) but with the condition $\sum_{i\in I} \gamma_i x_i < \infty$ replaced by $\sum_{i\in I} \gamma_i (x_i)^2 < \infty$.

The most basic facts about the process $(X(t))_{t\geq 0}$ are summarized in the following result.

THEOREM 1.1 (Shiga and Shimizu [42]).

(a) *The* SSDE *in* (1.2) *has a unique strong solution* $(X(t))_{t\geq 0}$ *on* $\mathcal{E}_1$ *with continuous paths.*

(b) $(X(t))_{t\geq 0}$ *is the unique Markov process on* $\mathcal{E}_1$ *whose semigroup* $(S(t))_{t\geq 0}$ *satisfies*

$$S(t)f - f = \int_0^t S(s)Lf\, ds, \qquad f \in C_0^2(\mathcal{E}_1), \tag{1.6}$$

*where $C_0^2(\mathcal{E}_1)$ is the space of functions on $\mathcal{E}_1$ depending on finitely many components and twice continuously differentiable in each component, and L is the pregenerator*

$$(Lf)(x) = \sum_{i\in I}\left\{\sum_{j\in I} a(i,j)[x_j - x_i]\right\}\frac{\partial f}{\partial x_i} + \frac{1}{2}\sum_{i\in I} b x_i^2 \frac{\partial^2 f}{\partial x_i^2}, \qquad x \in \mathcal{E}_1. \tag{1.7}$$

(c) *Restricted to $\mathcal{E}_2$, $(X(t))_{t\geq 0}$ is a diffusion process with the Feller property.*

The model defined by (1.2) represents a special case of the SSDE

$$dX_i(t) = \sum_{j\in I} a(i,j)[X_j(t) - X_i(t)]\,dt + \sqrt{g(X_i(t))}\,dW_i(t), \qquad i \in I, t \geq 0, \tag{1.8}$$

with $g$: $(-\infty,\infty) \to [0,\infty)$ some locally Lipschitz continuous function. This SSDE has, as far as its long-time behavior is concerned, four important classes:

(i) $g(x) > 0$ on $(0,1)$.
Examples: $g(x) = x(1-x)$ Fisher–Wright, $g(x) = (x(1-x))^2$ Ohta–Kimura.
(ii) $g(x) > 0$ on $(-\infty,\infty)$ and $g(x) = o(x^2)$ as $x \to \pm\infty$.
Example: $g(x) \equiv \sigma^2$ critical Ornstein–Uhlenbeck.
(iii) $g(x) > 0$ on $(0,\infty)$ and $g(x) = o(x^2)$ as $x \to \infty$.
Example: $g(x) = x$ Feller's continuous-state branching diffusion.
(iv) $g(x) > 0$ on $(0,\infty)$ and $g(x) \sim bx^2$ as $x \to \infty$.

Classes (i)–(iii) are well understood [9, 10, 19, 20, 24, 40, 41]. The qualitative properties of the process defined by (1.8) are similar for these three classes, and the universality of the long-time behavior as a function of $g$ has been systematically investigated via renormalization methods [1, 2, 13, 14, 15, 16, 26]. Class (iv), which is the subject of the current paper, is very different. For the case where $a(\cdot,\cdot)$ is simple random walk, this class was investigated in [41] and in the memoir by Carmona and Molchanov [7], where some of our results were already established.

The long-time behavior of the process defined by (1.2) is fairly complex. In order to keep the exposition transparent, we restrict our analysis to a subclass of models given by the following additional requirements:

$$\begin{aligned} &I = \mathbb{Z}^d, \qquad d \geq 1, \\ &a(\cdot,\cdot) \text{ is homogeneous: } a(i,j) = a(0,j-i) \qquad \forall i,j \in I, \\ &a(\cdot,\cdot) \text{ is irreducible: } \sum_{n=0}^{\infty}[a^n(i,j) + a^n(j,i)] > 0 \qquad \forall i,j \in I. \end{aligned} \tag{1.9}$$

Moreover, we put $a(0,0)=0$.

Before we start, let us fix some notation. We write $\mathcal{P}(\mathcal{E}_1)$ for the set of probability measures on $(\mathcal{E}_1, \mathcal{B}(\mathcal{E}_1))$, with $\mathcal{B}$ the Borel $\sigma$-algebra. For $\rho \in \mathcal{P}(\mathcal{E}_1)$, we write $E^\rho$ to denote expectation with respect to $\rho$. A measure $\rho \in \mathcal{P}(\mathcal{E}_1)$ is called *shift-invariant* if

$$\rho((X_i)_{i\in I} \in A) = \rho((X_{i+j})_{i\in I} \in A) \qquad \forall j \in I\ \forall A \in \mathcal{B}(\mathcal{E}_1), \tag{1.10}$$

is called *mixing* if

$$\lim_{\|k\|\to\infty} E^\rho(\, f\,[g \circ \sigma_k]\,) = E^\rho(f)\, E^\rho(g) \tag{1.11}$$

for all bounded $f, g : \mathcal{E}_1 \to \mathbb{R}$ that depend on finitely many coordinates, where $\sigma_k$ is the $k$-shift acting on $I$, and is called *associated* if

$$E^\rho(f_1 f_2) \geq E^\rho(f_1)\, E^\rho(f_2) \tag{1.12}$$

for all bounded $f_1, f_2 : \mathcal{E}_1 \to \mathbb{R}$ that depend on finitely many coordinates and that are nondecreasing in each coordinate.

We further need

$$\begin{aligned} \mathcal{T} &= \{\rho \in \mathcal{P}(\mathcal{E}_1) : \rho \text{ is shift-invariant}\}, \\ \mathcal{T}^1 &= \{\rho \in \mathcal{T} : E^\rho(X_0) < \infty\}, \end{aligned} \tag{1.13}$$

and

$$\mathcal{T}^1_\Theta = \{\rho \in \mathcal{T}^1 : \rho \text{ is shift ergodic}, E^\rho(X_0) = \Theta\}, \qquad \Theta \in [0,\infty). \tag{1.14}$$

The set of extreme points of a convex set $C$ is written $C_e$. The element $(x_i)_{i\in I}$ with $x_i = \Theta$ for all $i \in I$ is denoted by $\underline{\Theta}$. The initial distribution of our system is denoted by $\mu = \mathcal{L}(X(0))$ and is assumed to be concentrated on $\mathcal{E}_1$. The symbols $P, E$ without index denote probability and expectation with respect to $\mu$ and the Brownian motion driving (1.2). The notation w-lim means weak limit.

1.5. *Phase transitions.* In Theorems 1.2–1.5 below we state our main results on the long-time behavior of $(X(t))_{t\geq 0}$ and on the properties of its equilibria. Let

$$\mathcal{I} = \{\rho \in \mathcal{P}(\mathcal{E}_1) : \rho \text{ is invariant}\} \tag{1.15}$$

be the set of all equilibrium measures $\rho$ of (1.2), that is, $\rho S(t) = \rho$ for all $t \geq 0$. This set of course depends on $a(\cdot,\cdot)$ and $b$.

1.5.1. *Recurrent case.* The ergodic behavior of our system is simple when $\widehat{a}(\cdot,\cdot)$ defined by

$$\widehat{a}(i,j)=\tfrac{1}{2}[a(i,j)+a(j,i)], \qquad i,j\in I, \tag{1.16}$$

is recurrent. Namely, the process becomes extinct independently of the value of $b$.

THEOREM 1.2 (Shiga [41]). *If $\widehat{a}(\cdot,\cdot)$ is recurrent, then for every $b>0$ and every initial distribution $\mu\in\mathcal{T}^1$:*

$$\underset{t\to\infty}{\text{w-lim}}\,\mathcal{L}(X(t))=\delta_{\underline{0}}. \tag{1.17}$$

*Consequently, there exists no equilibrium in $\mathcal{T}^1$ other than $\delta_{\underline{0}}$, that is,*

$$\mathcal{I}\cap\mathcal{T}^1=\delta_{\underline{0}}. \tag{1.18}$$

Using the fact that if $\mu\in\mathcal{T}^1_\Theta$ then $E(X_i(t))=\Theta$ for all $i\in I$ and $t\geq 0$, we conclude from Theorem 1.2 that the system *clusters*, that is, on only few sites there is a nontrivial mass but at these sites the mass is very large (for $t$ large).

1.5.2. *Transient case*: *regimes* (I), (II) *and* (III). In the case where $\widehat{a}(\cdot,\cdot)$ is transient, the ergodic behavior of our system depends on the parameter $b$ and we observe interesting phase transitions. There are three regimes, separated by two critical values.

(I) *Small* $b$. Define the Green function

$$\widehat{G}(i,j)=\sum_{n=0}^{\infty}\widehat{a}^{\,n}(i,j), \qquad i,j\in I, \tag{1.19}$$

and put

$$b_2=\frac{2}{\widehat{G}(0,0)}. \tag{1.20}$$

We first consider the regime

$$\text{(I)} \qquad \widehat{a}(\cdot,\cdot)\ \text{transient}, \qquad b\in(0,b_2). \tag{1.21}$$

THEOREM 1.3 (Shiga [41], Cox and Greven [10]). *In regime* (I):

(a) *For $\mu=\delta_{\underline{\Theta}}$ with $\Theta\in[0,\infty)$ the following limit exists:*

$$\nu_\Theta=\underset{t\to\infty}{\text{w-lim}}\,\mathcal{L}(X(t)). \tag{1.22}$$

(b) *The measure* $\nu_\Theta$ *satisfies*

$$\begin{aligned}
&\nu_\Theta \in (\mathcal{I} \cap \mathcal{T}^1)_e,\\
&\nu_\Theta \text{ is shift-invariant, mixing and associated,}\\
&E^{\nu_\Theta}(X_0) = \Theta; \nu_\Theta \text{ is not a point mass if } \Theta > 0,\\
&E^{\nu_\Theta}([X_0]^2) < \infty.
\end{aligned} \tag{1.23}$$

(c) *The set of shift-invariant extremal equilibria is given by*

$$(\mathcal{I} \cap \mathcal{T}^1)_e = \{\nu_\Theta\}_{\Theta \in [0,\infty)}. \tag{1.24}$$

(d) *For every* $\mu \in \mathcal{T}^1_\Theta$ *with* $\Theta \in [0, \infty)$:

$$\underset{t\to\infty}{\text{w-lim}}\, \mathcal{L}(X(t)) = \nu_\Theta. \tag{1.25}$$

(e) *For every* $\mu \in \mathcal{T}_e$ *with* $E^\mu(X_0) = \infty$:

$$\underset{t\to\infty}{\text{w-lim}}\, \mathcal{L}(X(t)) = \delta_{\underline{\infty}}. \tag{1.26}$$

*Consequently,*

$$(\mathcal{I} \cap \mathcal{T})_e = (\mathcal{I} \cap \mathcal{T}^1)_e. \tag{1.27}$$

Theorem 1.3 tells us that if $b$ remains below an $a(\cdot,\cdot)$-dependent threshold, then the process $(X(t))_{t\geq 0}$ exhibits persistent behavior, in the sense that an equilibrium is approached with a spatial density equal to the initial spatial density and with a one-dimensional marginal that has a *finite* 2nd moment. This equilibrium is nontrivial unless the initial state is identically 0. If, on the other hand, the initial spatial density is infinite, then every component diverges in probability.

(II) *Moderate* $b$. We next consider the regime

$$\text{(II)} \qquad \widehat{a}(\cdot,\cdot) \text{ transient}, \qquad b \in [b_2, b_*). \tag{1.28}$$

In Section 4 we will obtain a variational expression for $b_*$ [see (4.19)]. This expression will turn out to be somewhat delicate to analyze.

THEOREM 1.4. *In regime* (II)*:*

(a) *The same properties hold as in Theorem* 1.3(a) *and* (c)–(e)*.*
(b) *The measure* $\nu_\Theta$ *satisfies*

$$\begin{aligned}
&\nu_\Theta \in (\mathcal{I} \cap \mathcal{T}^1)_e,\\
&\nu_\Theta \text{ is shift-invariant, mixing and associated,}\\
&E^{\nu_\Theta}(X_0) = \Theta; \nu_\Theta \text{ is not a point mass if } \Theta > 0,\\
&E^{\nu_\Theta}([X_0]^2) = \infty.
\end{aligned} \tag{1.29}$$

Theorem 1.4, which will be proved in Section 2, says that for moderate $b$ the equilibria continue to exist and to be well behaved, but with a one-dimensional marginal having *infinite* 2nd moment. The latter has important consequences for the fluctuations of the equilibrium in large blocks. Indeed, in regime (I) we may expect Gaussian limits after suitable scaling (see, e.g., Zähle [46, 47] in a different context), while in regime (II) we may expect non-Gaussian limits. In regime (II), the tail of $X_0$ under $\nu_\Theta$ is likely to be of stable law type, but a closer investigation of this question is beyond the scope of the present paper.

(III) *Large* $b$. Finally, we consider the regime

$$\text{(III)} \qquad \widehat{a}(\cdot,\cdot) \text{ transient}, \qquad b \in [b_*, \infty). \tag{1.30}$$

THEOREM 1.5. *In the interior of regime* (III), *for every* $\mu \in \mathcal{T}^1$:

$$\underset{t\to\infty}{\text{w-lim}}\, \mathcal{L}(X(t)) = \delta_{\underline{0}}. \tag{1.31}$$

*Consequently,* $\mathcal{I} \cap \mathcal{T}^1 = \delta_{\underline{0}}$.

We conjecture that there is local extinction at $b = b_*$.

Theorem 1.5, which will be proven in Section 4, shows that for large $b$ again *clustering* occurs, that is, the same situation as described in Theorem 1.2 for the case where $\widehat{a}(\cdot,\cdot)$ is recurrent.

1.6. *Finer analysis of the transient case.* In Section 1.5 we saw that different values of $b$ lead to qualitatively different behavior of the process $(X(t))_{t\geq 0}$. Therefore the question arises in which way the value of $b$ influences the properties of the process within one regime. For part of this finer analysis we need to assume that $a(\cdot,\cdot)$ is *symmetric*:

$$a(i,j) = a(j,i) \qquad \forall i,j \in I. \tag{1.32}$$

1.6.1. *Regime* (I). Let $\xi = (\xi(t))_{t\geq 0}$ be the random walk on $I$ with transition kernel $a(\cdot,\cdot)$ and jump rate 1, starting at 0. For $m \geq 2$, let $\xi^{(m)} = (\xi_1, \dots, \xi_m)$ be $m$ independent copies of $\xi$, and define the *differences random walk* $\eta^{(m)} = (\eta^{(m)}(t))_{t\geq 0}$ by putting

$$\eta^{(m)}(t) = (\xi_p(t) - \xi_q(t))_{1\leq p<q\leq m}. \tag{1.33}$$

This is a random walk on $I^{(m)}$, the subgroup of $I^{(1/2)m(m-1)}$ generated by all the possible pairwise differences of $m$ elements of $I$, with jump rate $m$ and transition kernel $a^{(m)}(\cdot,\cdot)$ that can be formally written out as

$$\begin{aligned} a^{(m)}(x,y) &= a^{(m)}(0, y-x) \\ &= \sum_{j\in I} a(0,j) \left[\frac{1}{m}\sum_{r=1}^{m} 1\{jD_r = y-x\}\right], \qquad x,y \in I^{(m)}, \end{aligned} \tag{1.34}$$

where $D_r$ is the triangular array of $-1, 0, +1$'s given by

$$D_r = (\delta_{pr} - \delta_{qr})_{1 \le p < q \le m} \tag{1.35}$$

and $jD_r$ denotes the triangular array obtained from $D_r$ by multiplying all its elements with the vector $j$. The factor $\frac{1}{m}$ comes from the fact that the $m$ random walks jump one at a time. Note that $a^{(m)}(\cdot,\cdot)$ is symmetric because of our assumption in (1.32). Note that $a^{(2)}(\cdot,\cdot) = \widehat{a}(\cdot,\cdot)$, the symmetrized transition kernel defined in (1.16), which is symmetric even without (1.32). The differences random walk is to be seen as the evolution of the random walks "relative to their center of mass." This will serve us later on.

Define the Green function

$$G^{(m)}(x,y) = \sum_{n=0}^{\infty} [a^{(m)}]^n(x,y), \qquad x,y \in I^{(m)}. \tag{1.36}$$

Also define the *collision function* $\sharp^{(m)} : I^{(m)} \to \mathbb{N}_0$ as

$$\begin{gathered} \sharp^{(m)}(z) = \sum_{1 \le p < q \le m} 1_{\{z_p - z_q = 0\}}, \\ z = (z_p - z_q)_{1 \le p < q \le m}, \qquad z_p, z_q \in I, \end{gathered} \tag{1.37}$$

and put

$$S^{(m)} = \mathrm{supp}(\sharp^{(m)}) \subset I^{(m)}. \tag{1.38}$$

Define

$$K^{(m)}(x,y) = \sqrt{\sharp^{(m)}(x)}\, G^{(m)}(x,y) \sqrt{\sharp^{(m)}(y)}, \qquad x,y \in S^{(m)}. \tag{1.39}$$

Viewed as an operator acting on $\ell^2(S^{(m)})$, $K^{(m)}$ is self-adjoint, positive and bounded. The latter two properties will be proved in Section 2.

The following result shows that in regime (I) there is an infinite sequence of critical values characterizing the convergence of successive moments.

THEOREM 1.6. *Suppose that $a(\cdot,\cdot)$ is symmetric. Then, in regime* (I), *there exists a sequence $b_2 \ge b_3 \ge b_4 \ge \cdots$ such that*:

(a) *If $\mu = \delta_{\underline{\Theta}}$ with $\Theta \in (0,\infty)$, then*

$$\lim_{t\to\infty} E([X_0(t)]^m) = E^{\nu_\Theta}([X_0]^m) \begin{cases} < \infty, & \text{for } b < b_m, \\ = \infty, & \text{for } b \ge b_m. \end{cases} \tag{1.40}$$

(b) *If $\mu = \delta_{\underline{\Theta}}$ with $\Theta \in (0,\infty)$, then*

$$\lim_{t\to\infty} \frac{1}{t} \log E([X_0(t)]^m) = \chi_m(b) \tag{1.41}$$

*exists with*

$$\chi_m(b) \begin{cases} = 0, & \text{for } b \le b_m, \\ > 0, & \text{for } b > b_m. \end{cases} \tag{1.42}$$

(c) *The critical value* $b_m$ *has the representation*

$$b_m = \frac{m}{\lambda_m} \tag{1.43}$$

*with* $\lambda_m \in (0,\infty)$ *the spectral radius of* $K^{(m)}$ *in* $\ell^2(S^{(m)})$*. This spectral radius is an eigenvalue if and only if* $b_{m-1} > b_m$*.*

(d) *The critical value* $b_2$ *is given by* (1.20)*, and*

$$b_2 \geq b_3 \geq b_4 \geq \cdots > 0. \tag{1.44}$$

*Moreover,*

$$\frac{2}{G^{(2)}(0,0)} = b_2 \leq (m-1)b_m \leq \frac{2}{G^{(m)}(0,0)} < 2, \tag{1.45}$$

*and* $\lim_{m\to\infty}(m-1)b_m$ *exists.*

(e) *The function* $b \mapsto \frac{1}{m}\chi_m(b)$ *is convex on* $[0,\infty)$ *and strictly increasing on* $[b_m,\infty)$*, with*

$$\lim_{b\to\infty} \frac{1}{bm}\chi_m(b) = \frac{1}{2}(m-1). \tag{1.46}$$

Theorem 1.6 will be proven in Section 3. Part (a) tells us that equilibria with finite $m$th moment exist if and only if $0 < b < b_m$. Part (b) tells us that the $m$th moment diverges exponentially fast if and only if $b > b_m$. The limit $\chi_m(b)$ is the *mth annealed Lyapunov exponent*. Part (c) gives a variational representation for $b_m$. Part (d) gives sharp bounds for $b_m$ and shows that the tail of the one-dimensional marginal of $\nu_\Theta$ decays algebraically with a power that is a nonincreasing function of $b$ when $b$ is small. It identifies the asymptotic behavior of this power as $\sim \mathrm{Cst}/b$ for $b \downarrow 0$. Part (e) shows that for large $b$ the curve $b \mapsto \frac{1}{m}\chi_m(b)$ has slope $\frac{1}{2}(m-1)$.

By Hölder's inequality, $m \mapsto \frac{1}{m}\chi_m(b)$ is nondecreasing. The system is called *intermittent of order* $n$ if

$$\frac{1}{n}\chi_n(b) < \frac{1}{n+1}\chi_{n+1}(b) < \frac{1}{n+2}\chi_{n+2}(b) < \cdots. \tag{1.47}$$

(It is shown in [7], Chapter III, that the first of these inequalities implies all the subsequent ones.) Thus, for all $n \geq 2$ our system is intermittent of order $n$ precisely when $b \in (b_{n+1}, b_n]$ (see also Figure 2 in Section 1.1).

We conjecture that $b_2 > b_3 > b_4 > \cdots$ [see open problem (A) in Section 1.2]. A partial result in this direction is the following:

COROLLARY 1.7. *Suppose that* $a(\cdot,\cdot)$ *is symmetric.*

(a) $(m-1)b_m \to 2$ *uniformly in* $m$ *as* $G^{(2)}(0,0) \to 1$.

(b) $b_2 > b_3 > \cdots > b_m$ *when* $G^{(2)}(0,0) \leq (m-1)/(m-2)$.

PROOF. (a) Obvious from (1.45).
(b) This follows from (1.45) and $G^{(m)}(0,0) > 1$. □

Claim (a) follows from (1.20) and (1.45), and corresponds to the limit when the random walk becomes more and more transient. This includes simple random walk on $\mathbb{Z}^d$ with $d \to \infty$. Thus, in this limit all inequalities in (1.44) become strict. Claim (b) follows from (1.45). This includes simple random walk on $\mathbb{Z}^d$ with $d \geq 3$.

As we will see in Sections 2–3, the representation for $b_m$ in (1.43) comes from a link with collision local time of random walks. Indeed, let

$$T(\xi^{(m)}) = \int_0^\infty \sum_{1 \leq p < q \leq m} 1_{\{\xi_p(t) = \xi_q(t)\}} \, dt \tag{1.48}$$

be the total collision local time (in pairs) of the $m$ independent copies of the random walk $\xi$. Then we will show that

$$b_m = \sup\{b > 0 : E^{\xi^{(m)}}(\exp[bT(\xi^{(m)})]) < \infty\}. \tag{1.49}$$

1.6.2. *Regime* (II). The next conjecture says that regime (II) is nonempty and may therefore be seen as an extension of Theorem 1.6(d).

CONJECTURE 1.8. *$b_* > b_2$ when $\widehat{a}(\cdot,\cdot)$ is transient.*

Conjecture 1.8 implies that *equilibria with stable law tails* occur in our system for moderate $b$ [recall (1.29)]. We conjecture that the system locally dies out at $b_*$ [see open problem (B) in Section 1.2].

In view of (1.49), we may ask whether it is possible to obtain a variational characterization for $b_*$. To that end, let $\xi = (\xi(t))_{t \geq 0}$ and $\xi' = (\xi'(t))_{t \geq 0}$ be two independent copies of the random walk on $I$ with transition kernel $a(\cdot,\cdot)$ and jump rate 1, both starting at 0. Let

$$T(\xi, \xi') = \int_0^\infty 1_{\{\xi(t) = \xi'(t)\}} \, dt \tag{1.50}$$

be their collision local time. Define

$$b_{**} = \sup\{b > 0 : E^{\xi'}(\exp[bT(\xi, \xi')]) < \infty \ \xi\text{-a.s.}\}, \tag{1.51}$$

where we note that $\{E^{\xi'}(\exp[bT(\xi, \xi')]) < \infty\}$ is a tail event for $\xi$. Since $b_2$ *is given by the same formula as* (1.51) *but with the average taken over both* $\xi$ *and* $\xi'$ [recall (1.49)], we have $b_{**} \geq b_2$. The proof of Conjecture 1.8 may be achieved by showing that

$$b_* \geq b_{**} \quad \text{and} \quad b_{**} > b_2. \tag{1.52}$$

In Section 4, we prove the first inequality and argue in favor of the second inequality. A full proof of the latter is deferred to [5].

1.6.3. *Regime* (III). The last result shows that in regime (III) the system gets extinct very rapidly.

THEOREM 1.9 (Cranston, Mountford and Shiga [12]). *In regime* (III):

(a) *If* $\mu = \delta_{\underline{\Theta}}$ *with* $\Theta \in (0, \infty)$, *then*

$$\lim_{t\to\infty} \frac{1}{t} \log X_0(t) = \chi_*(b) \tag{1.53}$$

*exists and is constant a.s.*

(b) *There exists a* $\tilde{b}_* \in [b_*, \infty)$ *such that*

$$\chi_*(b) \begin{cases} = 0, & \text{for } b \le \tilde{b}_*, \\ < 0, & \text{for } b > \tilde{b}_*. \end{cases} \tag{1.54}$$

(c)

$$\lim_{b\to\infty} \frac{\log b}{b} \left[\chi_*(b) + \frac{1}{2} b\right] \tag{1.55}$$

*exists in* $(0, \infty)$.

The limit $\chi_*(b)$ is the *quenched Lyapunov exponent*. Theorem 1.9 states that the speed of extinction is exponentially fast above a critical threshold $\tilde{b}_*$. Trivially,

$$\tilde{b}_* \ge b_*. \tag{1.56}$$

We conjecture that equality holds [see open problem (C) in Section 1.2]. See also Figures 1 and 2 in Section 1.1.

1.6.4. *Separation between regimes* (II) *and* (III). The key tool in the identification of $b_*$ is the notion of *Palm distribution* of our process $X$ at time $t$. This is the law of the process seen from a “randomly chosen mass” drawn at time $t$. This concept was introduced by Kallenberg [36] in the study of branching particle systems with migration. There the idea is to take a large box at time $t$, pick a particle at random from this box (the “tagged particle”), shift the origin to the location of this particle, consider the law of the shifted configuration, and let the box tend to infinity. Under suitable conditions, a limiting law is obtained, which is called the Palm distribution. Similarly, in our system the Palm distribution is a size biasing of the original distribution according to the “mass” at the origin. The criterion for *survival versus extinction of the original distribution translates into tightness versus divergence of the Palm distribution.*

This criterion is useful for two reasons. First, the size biasing is an easy operation. Second, often it is possible to obtain a representation formula for the Palm distribution in terms of a nice Markov process. For instance, for branching particle systems the Palm distribution is obtained from an independent superposition of the original distribution and a realization of the so-called Palm canonical distribution. The latter can be identified as a branching random walk with immigration of particles at rate 1 along the path of the tagged particle. Fortunately, we can give an explicit representation of the Palm distribution of our process $X$ as well, namely, as the solution of a system of biased stochastic differential equations (see Section 4 for details). It turns out that the latter again has a (Feynman–Kac type) representation formula for the single components as an expectation over an exponential functional of the Brownian motions, the random walk, and an additional tagged random walk, with the expectation running over the two random walks.

We will use the Palm distribution to identify $b_*$. We will see that, within the interval $(b_2, b_*)$, we can distinguish between a regime where the average of the Palm distribution over the Brownian motions (i.e., the Palm distribution conditioned on the tagged path) is tight as $t \to \infty$ and a regime where it diverges. The separation between these two regimes is $b_{**}$. Within the interval $[b_{**}, b_*)$, we can separate further by conditioning the Palm distribution also on the Brownian increments along the tagged path. However, we will not pursue this point further, even though it is of interest for a better insight into what controls our system. See Birkner [3, 4] for more background.

We will see in Section 4.1.2 that (1.50) plays an important role in the description of the Palm distribution. Equation (4.19) in Section 4.1.2 identifies $b_*$. However, this formula is much harder than the one for $b_{**}$ in (1.51). It would be interesting to know whether there exists a characterization of $b_*$ in terms of the collision local time of random walks [see open problem (D) in Section 1.2], in the same way as for $b_m$ in (1.49) and for $b_{**}$ in (1.51).

## 2. Moment calculations.

2.1. *Definition of* $b_m$, $\bar{b}_m$ *and* $\widetilde{b}_m$. For $\Theta > 0$ and $m \geq 2$, let

$$
\begin{aligned}
b_m &= \sup\{b > 0 : \nu_{\underline{\Theta}} \neq \delta_{\underline{0}}, E^{\nu_\Theta}([X_0]^m) < \infty\},\\
\bar{b}_m &= \sup\Big\{b > 0 : \limsup_{t\to\infty} E([X_0(t)]^m \mid X(0) = \underline{\Theta}) < \infty\Big\}, \qquad (2.1)\\
\widetilde{b}_m &= \sup\Big\{b > 0 : \limsup_{t\to\infty} \frac{1}{t} \log E([X_0(t)]^m \mid X(0) = \underline{\Theta}) = 0\Big\}.
\end{aligned}
$$

In these definitions the choice of $\Theta$ is irrelevant as long as $\Theta > 0$, as is evident from Lemma 2.1 below.

In Section 2.2 we derive a representation formula for the solution of (1.2), which is due to Shiga [41] and which plays a key role in the present paper. We also derive a self-duality property, which is needed to obtain convergence to equilibrium. In Section 2.3 we express the $m$th moment of a single component of our system, at time $t$, in terms of the collision local time, up to time $t$, of $m$ independent copies of our random walk. In Section 2.4 we prove that $\nu_\Theta$ exists and that $b_m = \bar{b}_m$. In Section 2.5 we prove some basic properties of $G^{(m)}$ and $K^{(m)}$ defined in (1.36) and (1.39). The results in this section will be used in Sections 3–4 to prove Theorems 1.4–1.5, 1.6 and 1.9.

2.2. *Representation formula and self-duality.* If our process starts in a constant initial configuration, then a nice (Feynman–Kac type) representation formula is available. This formula will play a key role throughout the paper.

LEMMA 2.1. *The process $(X(t))_{t\geq 0}$ starting in $X(0) = \underline{\Theta}$ can be represented as the following functional of the Brownian motions:*

$$X_i(t) = \Theta e^{-(1/2)bt} E_i^\xi\left(\exp\left[\sqrt{b}\int_0^t \sum_{j\in I} 1_{\{\xi(t-s)=j\}}\, dW_j(s)\right]\right), \qquad i \in I, t \geq 0, \tag{2.2}$$

*where $\xi = (\xi(t))_{t\geq 0}$ is the random walk on $I$ with transition kernel $a(\cdot,\cdot)$ and jump rate 1, and the expectation is over $\xi$ conditioned on $\xi(0) = i$ ($\xi$ and $W$ are independent).*

PROOF. This lemma appears in [41] without proof. We write out the proof here, because it will serve us later on. The symbol $\mathbb{1}$ denotes the identity matrix. Note that $\Theta$ enters into (2.2) only as a front factor.

FACT 1. *For all $i$ and $t$:*

$$X_i(t) = \Theta + \sqrt{b}\int_0^t \sum_j a_{t-s}(i,j)\, X_j(s)\, dW_j(s) \tag{2.3}$$

*with $a_t = \exp[t(a - \mathbb{1})]$.*

PROOF. Fix $i$ and $t$. For $0 \leq s \leq t$, let $Y_i(s) = \sum_j a_{t-s}(i,j) X_j(s)$. [The infinite sum is finite due to the fact that, by Theorem 1.1(a), $(X(t))_{t\geq 0}$ lives in $\mathcal{E}_1$ defined by (1.4).] Then

$$dY_i(s) = \left[\sum_j (\mathbb{1} - a) a_{t-s}(i,j)\, X_j(s)\right] ds + \sum_j a_{t-s}(i,j)\, dX_j(s). \tag{2.4}$$

From (1.2) we have

$$dX_i(s) = \sum_j (a - \mathbb{1})(i,j) X_j(s)\, ds + \sqrt{b}\, X_i(s)\, dW_i(s), \tag{2.5}$$

which after substitution into (2.4) and cancellation of two terms gives

$$dY_i(s) = \sqrt{b} \sum_j a_{t-s}(i,j) X_j(s)\, dW_j(s). \tag{2.6}$$

Integrate both sides from 0 to $t$, and note that $Y_i(0) = \sum_j a_t(i,j)X_j(0) = \Theta$ and $Y_i(t) = \sum_j a_0(i,j)X_j(0) = X_i(t)$, to get the claim. □

FACT 2. *For all $t$:*

$$\exp[\sqrt{b} Z_t(t) - \tfrac{1}{2}bt] = 1 + \sqrt{b} \int_0^t \exp[\sqrt{b}\, Z_t(s) - \tfrac{1}{2}bs]\, dZ_t(s) \tag{2.7}$$

*with* $Z_t(s) = \int_0^s \sum_j 1_{\{\xi(t-r)=j\}}\, dW_j(r)$.

PROOF. Fix $t$. For $z \in \mathbb{R}$ and $s \geq 0$, let $f(z,s) = e^{\sqrt{b}\,z-(1/2)bs}$ and put $g(s) = f(Z_t(s), s)$. Itô's formula gives

$$\begin{aligned} dg(s) &= f_z(Z_t(s), s)\, dZ_t(s) + \tfrac{1}{2} f_{zz}(Z_t(s), s)(dZ_t(s))^2 \\ &\quad + f_s(Z_t(s), s)\, ds, \end{aligned} \tag{2.8}$$

which after cancellation of two terms [because $\frac{1}{2}f_{zz} + f_s = 0$ and $(dZ_t(s))^2 = ds$] gives

$$dg(s) = g(s)\sqrt{b}\; dZ_t(s). \tag{2.9}$$

Integrate both sides from 0 to $t$ and use that $g(0) = f(Z_t(0), 0) = 1$, to get the claim. □

The proof of the representation formula in Lemma 2.1 is now completed as follows. Let $\widetilde{X}_i(t)$ denote the right-hand side of (2.2). Taking the expectation over $\xi$ conditioned on $\xi(0) = i$ on both sides of (2.7), we get

$$\begin{aligned} \Theta^{-1}\widetilde{X}_i(t) &= 1 + \sqrt{b} \int_0^t E_i^{\xi}\left( \sum_j 1_{\{\xi(t-s)=j\}} \Theta^{-1}\widetilde{X}_j(s)\, dW_j(s) \right) \\ &= 1 + \sqrt{b} \int_0^t \sum_j a_{t-s}(i,j)\; \Theta^{-1}\widetilde{X}_j(s)\, dW_j(s), \end{aligned} \tag{2.10}$$

where the first equality uses the Markov property of $\xi$ at time $t - s$. Thus we see that $\widetilde{X}_i(t)$ satisfies (2.3). Since $\widetilde{X}_i(0) = \Theta = X_i(0)$ for all $i$, we may therefore conclude that $\widetilde{X}_i(t) = X_i(t)$ for all $i$ and $t$, by the strong uniqueness of the solution of our system (1.2) [recall Theorem 1.1(a)]. □

In addition to the representation formula in Lemma 2.1, we have another nice property: our process is self-dual. Let

$$a^*(i,j) = a(j,i), \qquad i,j \in I, \tag{2.11}$$

be the reflected transition kernel. Let (1.2*) denote (1.2) with $a(\cdot,\cdot)$ replaced by $a^*(\cdot,\cdot)$. Abbreviate $\langle x, x^*\rangle = \sum_{i\in I} x_i x_i^*$.

LEMMA 2.2. *Let $X = (X(t))_{t\geq 0}$ be the solution of* (1.2) *starting from any $X(0) \in \mathcal{E}_1$. Let $X^* = (X^*(t))_{t\geq 0}$ be the solution of* (1.2*) *starting from any $X^*(0) \in \mathcal{E}_1$ such that $\langle \underline{1}, X^*(0)\rangle < \infty$. Then*

$$E^X(e^{-\langle X(t), X^*(0)\rangle}) = E^{X^*}(e^{-\langle X(0), X^*(t)\rangle}) \qquad \forall t \geq 0. \tag{2.12}$$

PROOF. See Cox, Klenke and Perkins [11]. □

2.3. *Representation of the mth moment in terms of collision local time.* Let us abbreviate

$$W = (\{W_i(t)\}_{i\in I})_{t\geq 0} \tag{2.13}$$

and write

$$E([X_0(t)]^m \mid X(0) = \underline{\Theta}) = E^W([X_0(t)]^m \mid X(0) = \underline{\Theta}) \tag{2.14}$$

to display that (1.2) is driven by $W$. This subsection contains a moment calculation in which we use the representation formula of Lemma 2.1 to express the right-hand side of (2.14) as the expectation of the exponential of $b$ times the *collision local time* of $m$ independent copies of the random walk with transition kernel $a(\cdot,\cdot)$ and jump rate 1, all starting at 0.

We begin by checking that the evolution is mean-preserving. This property is evident from (2.3), but its proof will serve as a preparation for the calculation of the higher moments.

LEMMA 2.3. $E^W(X_0(t) \mid X(0) = \underline{\Theta}) = \Theta$ *for all* $t \geq 0$.

PROOF. Taking the expectation over $W$ in (2.2) and using Fubini's theorem, we have

$$\begin{aligned} &E^W(X_0(t) \mid X(0) = \Theta) \\ &\quad = \Theta e^{-(1/2)bt} E_0^\xi\left(E^W\left(\exp\left[\sqrt{b}\int_0^t \sum_{i\in I} 1_{\{\xi(t-s)=i\}}\, dW_i(s)\right]\right)\right). \end{aligned} \tag{2.15}$$

Since the Brownian motions $W$ are i.i.d. and have independent increments, it follows that for any $\xi$:

$$\int_0^t \sum_{i\in I} 1_{\{\xi(t-s)=i\}}\, dW_i(s) \triangleq W'(t), \tag{2.16}$$

where $W' = (W'(t))_{t\geq 0}$ is a single Brownian motion and $\triangleq$ denotes equality in distribution. Combining (2.15) and (2.16) we arrive at (the expectation over $\xi$ being irrelevant)

$$E^W(X_0(t) \mid X(0) = \underline{\Theta}) = \Theta e^{-(1/2)bt} E^{W'}(\exp[\sqrt{b}W'(t)]). \tag{2.17}$$

Now use that, by Itô's formula, $\exp[\sqrt{b}W'(t) - \frac{1}{2}bt]$ is a martingale, to get that the r.h.s. of (2.17) equals $\Theta$. □

A version of the above argument will produce the following expression for the moments of order $m \geq 2$.

Lemma 2.4. *Let $\xi^{(m)} = (\xi_1, \dots, \xi_m)$ be $m$ independent copies of the random walk with transition kernel $a(\cdot,\cdot)$ and jump rate 1, all starting at 0. Then*

$$E^W([X_0(t)]^m \mid X(0) = \underline{\Theta}) = \Theta^m E^{\xi^{(m)}}(\exp[bT^{(m)}(t)]), \tag{2.18}$$

*where*

$$\begin{aligned} T^{(m)}(t) &= \sum_{1\leq k<l\leq m} T_{kl}(t), \\ T_{kl}(t) &= \int_0^t ds 1_{\{\xi_k(s)=\xi_l(s)\}}\, ds, \end{aligned} \tag{2.19}$$

*is the collision local time (in pairs) up to time $t$.*

Proof. Similarly as in (2.15) we may use (2.2) to write

$$\begin{aligned} &E^W([X_0(t)]^m \mid X(0) = \underline{\Theta}) \\ &\quad = \Theta^m e^{-(m/2)bt} E^{\xi^{(m)}}\left(E^W\left(\exp\left[\sqrt{b}\int_0^t \sum_{k=1}^m \sum_{i\in I} 1_{\{\xi_k(t-s)=i\}}\, dW_i(s)\right]\right)\right). \end{aligned} \tag{2.20}$$

Next, let $W^{(m)} = (W'_1, \dots, W'_m)$ be $m$ independent Brownian motions. Then the analogue of (2.16) reads

$$\begin{aligned} &\int_0^t \sum_{k=1}^m \sum_{i\in I} 1_{\{\xi_k(t-s)=i\}}\, dW_i(s) \\ &\quad \triangleq \int_0^t \sum_{p=1}^m \sum_{j^{(p)}} 1_{j^{(p)}}(t-s) \sum_{q=1}^p j_q\, dW'_q(s). \end{aligned} \tag{2.21}$$

Here $1_{j^{(p)}}(t-s)$ denotes the indicator of the event that at time $t-s$ the components of $\xi^{(m)} = (\xi_1, \dots, \xi_m)$ coincide in $p$ subgroups of sizes $j^{(p)} = (j_1, \dots, j_p)$ with $j_1 + \cdots + j_p = m$. The equality in (2.21) again follows from the fact that the Brownian motions $W$ are i.i.d. and have independent increments. The point to note here is that all $j_q$ random walks in the $q$th coincidence group pick up the same increment of the Brownian motion in $W$ at the site where they coincide at time $s$, and this increment has the same distribution as $dW'_q(s)$. Next, define

$$T_{j^{(p)}}(t) = \int_0^t 1_{j^{(p)}}(t-s)\,ds. \tag{2.22}$$

Then clearly we have

$$\int_0^t 1_{j^{(p)}}(t-s)\,dW'_q(s) \triangleq W'_q(T_{j^{(p)}}(t)). \tag{2.23}$$

Now combine (2.20)–(2.23) to get

$$\begin{aligned} &E^W([X_0(t)]^m \mid X(0) = \underline{\Theta}) \\ &\quad = \Theta^m e^{-(m/2)bt} E^{\xi^{(m)}}\left(E^{W^{(m)}}\left(\exp\left[\sqrt{b}\sum_{p=1}^m \sum_{j^{(p)}} \sum_{q=1}^p j_q W'_q(T_{j^{(p)}}(t))\right]\right)\right) \qquad (2.24)\\ &\quad = \Theta^m e^{-(m/2)bt} E^{\xi^{(m)}}\left(\exp\left[b\sum_{p=1}^m \sum_{j^{(p)}} T_{j^{(p)}}(t)\left\{\sum_{q=1}^p \tfrac{1}{2} j_q^2\right\}\right]\right). \end{aligned}$$

Finally, absorb the term $-\frac{m}{2}bt$ into the sum by writing $\frac{1}{2}j_q(j_q-1)$ instead of $\frac{1}{2}j_q^2$ [use that $\sum_{p=1}^m \sum_{j^{(p)}} T_{j^{(p)}}(t) = t$ and $\sum_{q=1}^p j_q = m$]. The resulting exponent is the same as $b$ times the collision local time in (2.19). □

2.4. *Convergence to equilibrium and* $b_m = \bar{b}_m$. The following important facts will be needed later on and will be derived via the representation formula in Lemma 2.1 and the self-duality in Lemma 2.2.

PROPOSITION 2.5. *For all* $\widehat{a}(\cdot,\cdot)$ *transient and all* $b > 0$:

(a) $\nu_\Theta = \text{w-lim}_{t\to\infty} \mathcal{L}(X(t) \mid X(0) = \underline{\Theta})$ *exists, is shift-invariant and associated for all* $\Theta \in [0,\infty)$.

(b) $\nu_\Theta$ *is mixing for all* $\Theta \in [0,\infty)$.

(c) $b_m = \bar{b}_m$ *for all* $m \geq 2$.

(d) $\lim_{t\to\infty} E([X_0(t)]^m \mid X(0) = \underline{\Theta}) = E^{\nu_\Theta}([X_0]^m)$ *for all* $\Theta \in [0,\infty)$ *and all* $m \geq 2$.

PROOF. (a) The proof of existence uses Lemma 2.2. If $X(0) = \underline{\Theta}$ and $X^*(0) = f$, then (2.12) reads

$$\begin{aligned} &E^X(e^{-\langle X(t), f\rangle} \mid X(0) = \underline{\Theta}) \\ &\qquad = E^{X^*}(e^{-\Theta M(t)} \mid M(0) = \langle \underline{1}, f\rangle) \end{aligned} \tag{2.25}$$

with

$$M(t) = \langle \underline{1}, X^*(t) \rangle. \tag{2.26}$$

Since $(M(t))_{t\geq 0}$ is a nonnegative martingale [as is obvious from (2.3) with the reflected transition kernel], we have that $\lim_{t\to\infty} M(t) = M(\infty) < \infty$ $W$-a.s. for every $f$ [use that $E^{X^*}(M(\infty)) \leq M(0) = \langle \underline{1}, f \rangle < \infty$]. Hence we conclude that $X(t)$ converges in law to a limit, which we call $\nu_{\underline{\Theta}}$, given by

$$\begin{aligned} \int e^{-\langle x,f \rangle} \nu_{\underline{\Theta}}(dx) &= E^{X^*}(e^{-\Theta M_\infty} \mid M(0) = \langle \underline{1}, f \rangle) \\ &\qquad \text{for all } f \text{ such that } \langle \underline{1}, f \rangle < \infty. \end{aligned} \tag{2.27}$$

Because $\delta_{\underline{\Theta}}$ is shift-invariant, so is $\nu_\Theta$. The fact that $\nu_\Theta$ is associated follows from Cox and Greven [10]. There it is shown that for systems of the type in (1.2)—even with a general diffusion term—the evolution preserves the associatedness. Since $\delta_{\underline{\Theta}}$ is associated, the system is associated at time zero and hence at all later times, and the equilibrium inherits this property.

(b) For $0 < b < b_2$ the mixing property of the equilibrium $\nu_\Theta$ was proved in Cox and Greven [10] via a covariance argument. However, for $b_2 \leq b < b_*$ covariances are infinite, and so we must follow a different route.

The proof uses the exponential duality in Lemma 2.2. We will prove that, for all $f, g \in \mathcal{E}_1$ [recall (1.4)] with finite support,

$$\lim_{\|k\|\to\infty} E^{\nu_\Theta}(e^{-\langle X,f \rangle} e^{-\langle X,\sigma_k g \rangle}) = E^{\nu_\Theta}(e^{-\langle X,f \rangle}) \, E^{\nu_\Theta}(e^{-\langle X,g \rangle}), \tag{2.28}$$

where $\sigma_k g = g \circ \sigma_k$ with $\sigma_k$ the $k$-shift acting on $I$. This implies the mixing property, because the Laplace functional determines the distribution.

*Step* 1. In order to prove (2.28), we use the self-duality of our process and the fact that $\nu_\Theta = \text{w-}\lim_{t\to\infty} \delta_{\underline{\Theta}} S(t)$, as follows. Denote by

$$X^{*,h} = (X^{*,h}(t))_{t\geq 0} = (\{X_i^{*,h}(t)\}_{i\in\mathbb{Z}^d})_{t\geq 0} \tag{2.29}$$

our process with reflected transition kernel $a^*(\cdot,\cdot)$ [recall (2.11)] starting from initial configuration $h \in \mathcal{E}_1$ with finite support. Then

$$\begin{aligned} E^{\nu_\Theta}(e^{-\langle X,f \rangle} e^{-\langle X,g \rangle}) &= E^{\nu_\Theta}(e^{-\langle X,f+g \rangle}) \\ &= \lim_{t\to\infty} E(e^{-\langle X(t),f+g \rangle} \mid X(0) = \underline{\Theta}) \\ &= \lim_{t\to\infty} E(e^{-\langle \underline{\Theta}, X^{*,f+g}(t) \rangle}). \end{aligned} \tag{2.30}$$

Observe that, by the linearity of the system, we may use the same Brownian motions for $X^{*,f}$ and $X^{*,g}$, which gives us in addition

$$X^{*,f+g} \triangleq X^{*,f} + X^{*,g}. \tag{2.31}$$

Hence, in order to verify (2.28), we must investigate the quantity

$$\lim_{t\to\infty} E(e^{-\langle \underline{\Theta}, X^{*,f}(t)+X^{*,\sigma_k g}(t)\rangle}) \tag{2.32}$$

and show that it factorizes in the limit as $\|k\|\to\infty$.

Next, note that

$$(\langle \underline{\Theta}, X^{*,f}(t)\rangle)_{t\geq 0} \quad \text{and} \quad (\langle \underline{\Theta}, X^{*,\sigma_k g}(t)\rangle)_{t\geq 0} \tag{2.33}$$

are (continuous-path square-integrable) nonnegative martingales. In particular, their limit as $t\to\infty$ exists by the martingale convergence theorem. Their covariation over the time interval $[0,\infty)$ is given by

$$\begin{aligned} C(f,\sigma_k g) &= \int_0^\infty ds \langle X^{*,f}(s), X^{*,\sigma_k g}(s)\rangle \\ &= \int_0^\infty ds \sum_{i\in I} X_i^{*,f}(s) X_i^{*,\sigma_k g}(s). \end{aligned} \tag{2.34}$$

Due to these structural properties, we know that if

$$C(f,\sigma_k g)\to 0 \qquad \text{in probability as } \|k\|\to\infty, \tag{2.35}$$

then the two martingales in (2.33) become independent as $\|k\|\to\infty$. Consequently, the two random variables

$$\lim_{t\to\infty} \langle \underline{\Theta}, X^{*,f}(t)\rangle \quad \text{and} \quad \lim_{t\to\infty} \langle \underline{\Theta}, X^{*,\sigma_k g}(t)\rangle \tag{2.36}$$

also become independent as $\|k\|\to\infty$, which proves (2.28) via (2.30) and (2.31). In order to prove (2.35), observe that, by the linearity of the system, it suffices to verify (2.28) for the special case where $f$ and $g$ are indicators of a single site in $I$, say, $p$ and $q$, respectively.

*Step* 2. Let $\xi$ and $\xi'$ be two independent random walks with transition kernel $a^*(\cdot,\cdot)$ and jump rate 1, both starting in $i\in I$. Then, for $f=1_{\{p\}}$ and $g=1_{\{q\}}$, it follows from Lemma 2.1 that

$$\begin{aligned} &X_i^{*,f}(s) X_i^{*,\sigma_k g}(s) \\ &\quad = e^{-bs} E_{i,i}^{\xi,\xi'}\Bigg(1\{\xi(s)=p,\xi'(s)=q+k\} \\ &\qquad\qquad \times \exp\Bigg[\sqrt{b}\int_0^s du \sum_{m\in I}[1\{\xi(s-u)=m\} \\ &\qquad\qquad\qquad + 1\{\xi'(s-u)=m\}]\, dW_m(u)\Bigg]\Bigg). \end{aligned} \tag{2.37}$$

Reversing time, we may start $\xi$ in $p$ and $\xi'$ in $q+k$, and give them transition kernel $a(\cdot,\cdot)$ and jump rate 1. Then

$$
\begin{aligned}
&X_i^{*,f}(s)X_i^{*,\sigma_k g}(s)\\
&\quad = e^{-bs}E^{\xi,\xi'}_{p,q+k}\biggl(1\{\xi(s)=\xi'(s)=i\}\\
&\qquad\qquad \times \exp\biggl[\sqrt{b}\int_0^s du \sum_{m\in I}[1\{\xi(u)=m\}\\
&\qquad\qquad\qquad + 1\{\xi'(u)=m\}]\,dW_m(u)\biggr]\biggr)
\end{aligned}
\tag{2.38}
$$

and so

$$
\begin{aligned}
&\langle X^{*,f}(s), X^{*,\sigma_k g}(s)\rangle\\
&\quad = \sum_{i\in I} X_i^{*,f}(s)X_i^{*,\sigma_k g}(s)\\
&\quad = e^{-bs}E^{\xi,\xi'}_{p,q+k}\biggl(1\{\xi(s)=\xi'(s)\}\\
&\qquad\qquad \times \exp\biggl[\sqrt{b}\int_0^s du \sum_{m\in I}[1\{\xi(u)=m\}\\
&\qquad\qquad\qquad + 1\{\xi'(u)=m\}]\,dW_m(u)\biggr]\biggr).
\end{aligned}
\tag{2.39}
$$

Next, for any $s$ we have

$$
\begin{aligned}
\langle X^{*,f}(s), X^{*,\sigma_k g}(s)\rangle &= \sum_{i\in I} X_i^{*,f}(s)X_i^{*,\sigma_k g}(s)\\
&\le \tfrac12\sum_{i\in I}[X_i^{*,f}(s)]^2 + \tfrac12\sum_{i\in I}[X_i^{*,\sigma_k g}(s)]^2\\
&= \tfrac12 M_1(s) + \tfrac12 M_2(s),
\end{aligned}
\tag{2.40}
$$

where

$$
\begin{aligned}
M_1(s) &= e^{-bs}E^{\xi}_{p}\biggl(\exp\biggl[2\sqrt{b}\int_0^s du \sum_{m\in I} 1\{\xi(u)=m\}\,dW_m(u)\biggr]\biggr),\\
M_2(s) &= e^{-bs}E^{\xi}_{q+k}\biggl(\exp\biggl[2\sqrt{b}\int_0^s du \sum_{m\in I} 1\{\xi(u)=m\}\,dW_m(u)\biggr]\biggr).
\end{aligned}
\tag{2.41}
$$

Below we will show that

$$\int_0^\infty ds\, M_i(s) < \infty \qquad W\text{-a.s. for } i = 1, 2. \tag{2.42}$$

Assuming (2.42), we pick $T > 0$ and estimate, with the help of (2.40),

$$\begin{aligned} C(f, \sigma_k g) &= \int_0^\infty ds \langle X^{*,f}(s), X^{*,\sigma_k g}(s)\rangle \\ &\leq \int_0^T ds\, \langle X^{*,f}(s), X^{*,\sigma_k g}(s)\rangle + \tfrac{1}{2}\int_T^\infty [M_1(s) + M_2(s)]. \end{aligned} \tag{2.43}$$

By (2.42) and the fact that the law of $(M_2(s))_{s\geq 0}$ is independent of $k$, it now suffices to show that

$$\int_0^T ds \langle X^{*,f}(s), X^{*,\sigma_k g}(s)\rangle \to 0 \qquad \text{in probability as } \|k\| \to \infty \text{ for any } T > 0. \tag{2.44}$$

*Step* 3. To prove (2.44), we return to (2.39). Write

$$\begin{aligned} &\sum_{m\in I} [1\{\xi(u) = m\} + 1\{\xi'(u) = m\}]\, dW_m(u) \\ &\quad = \sum_{m\in I} 1\{\xi(u) = \xi'(u) = m\} 2\, dW_m(u) \\ &\qquad + \sum_{m\in I} 1\{\xi(u) = m \neq \xi'(u)\}\, dW_m(u) \\ &\qquad + \sum_{m\in I} 1\{\xi(u) \neq m = \xi'(u)\}\, dW_m(u). \end{aligned} \tag{2.45}$$

Since the three terms in the right-hand side of (2.45) involve disjoint time intervals and the $W_m$'s have independent increments, it follows from (2.45) that, for any $\xi$,

$$\begin{aligned} &\sqrt{b}\int_0^s du \sum_{m\in I} [1\{\xi(u) = m\} + 1\{\xi'(u) = m\}]\, dW_m(u) \\ &\quad \triangleq W_1(4bT_s(\xi, \xi')) + W_2(b[s - T_s(\xi, \xi')]) + W_3(b[s - T_s(\xi, \xi')]), \end{aligned} \tag{2.46}$$

where $W_1, W_2, W_3$ are three independent Brownian motions, and $T_s(\xi, \xi') = \int_0^s du\, 1\{\xi(u) = \xi'(u)\}$ is the collision local time of $\xi$ and $\xi'$ up to time $s$. By combining (2.39) and (2.46), taking the expectation over $W$ (i.e., over $W_1, W_2, W_3$) and using Fubini's theorem, we get

$$E^W\left(\int_0^T ds \langle X^{*,f}(s), X^{*,\sigma_k g}(s)\rangle\right)$$

$$
\begin{aligned}
&= \int_0^T E^{\xi,\xi'}_{p,q+k}(1\{\xi(s)=\xi'(s)\}e^{bT_s(\xi,\xi')}) \\
&\le e^{bT} E^{\xi,\xi'}_{p,q+k}(T_T(\xi,\xi')).
\end{aligned} \tag{2.47}
$$

Clearly, for fixed $T$ the right-hand side tends to zero as $\|k\| \to \infty$, because $T_T(\xi,\xi') \le T$ and $T_T(\xi,\xi') \to 0$ in probability with respect to $\xi,\xi'$ as $\|k\| \to \infty$ for any fixed $T$.

*Step* 4. It remains to prove (2.42), which goes as follows. Let

$$
(M(s))_{s\ge 0} \qquad \text{with } M(s) = \langle \underline{1}, X^{*,f}(s)\rangle. \tag{2.48}
$$

This is a (continuous-path square-integrable) nonnegative martingale starting from a strictly positive and finite value (because $f = 1_{\{p\}}$ has finite support). From the dual of (1.2) (recall Lemma 2.2), we have

$$
\begin{aligned}
dM(s) &= \sum_{i\in I} dX^{*,f}_i(s) \\
&= \sum_{i,j\in I} a^*(i,j)[X^{*,f}_j(s) - X^{*,f}_i(s)]\,ds \\
&\quad + \sum_{i\in I} \sqrt{b[X^{*,f}_i(s)]^2}\,dW_i(s).
\end{aligned} \tag{2.49}
$$

The first term in the right-hand side is zero because $a^*(\cdot,\cdot)$ is doubly stochastic (being a random walk transition kernel [recall (1.9)] and $\langle \underline{1}, X^{*,f}(s)\rangle < \infty$. Hence

$$
M(s) \triangleq \widehat{W}(\tau(s)) \tag{2.50}
$$

with

$$
\tau(s) = \int_0^s du\, b \sum_{i\in I} [X^{*,f}_i(u)]^2 \tag{2.51}
$$

and $\widehat{W}$ a Brownian motion adapted to the filtration of $X^{*,f}$. By the martingale convergence theorem, we have

$$
\lim_{s\to\infty} M(s) = M(\infty) < \infty, \qquad W\text{-a.s.} \tag{2.52}
$$

[use that $E^{X^{*,f}}(M(\infty)) \le M(0) = \langle \underline{1}, f\rangle < \infty$]. Combining this with (2.50), we conclude that

$$
\lim_{s\to\infty} \tau(s) = \tau(\infty) < \infty, \qquad W\text{-a.s.} \tag{2.53}
$$

This completes the proof of (2.42), hence of (2.35), and therefore also of (2.28). □

(c) Fatou's lemma in combination with part (a) shows that

$$\liminf_{t\to\infty} E([X_0(t)]^m \mid X(0) = \underline{\Theta}) \geq E^{\nu_\Theta}([X_0]^m). \tag{2.54}$$

Hence $\bar{b}_m \leq b_m$. The converse is proved as follows. Assume that $\nu_\Theta \neq \delta_{\underline{0}}$. Define the $m$-point correlation function in equilibrium,

$$f(j_1, \dots, j_m) = E^{\nu_\Theta}\left(\prod_{p=1}^{m} X_{j_p}\right), \qquad j_1, \dots, j_m \in I, \tag{2.55}$$

where the indices need not be distinct. Since $\nu_\Theta$ is associated and shift-invariant [which was proved in part (a)], we have

$$\Theta^m \leq f(j_1, \dots, j_m) \leq f(0, \dots, 0) = E^{\nu_\Theta}([X_0]^m). \tag{2.56}$$

Moreover, from the equilibrium property of $\nu_\Theta$ we deduce that, for any $t > 0$,

$$\begin{aligned} E^{\nu_\Theta}\left(\prod_{p=1}^{m} X_{j_p}\right) &= \int \nu_\Theta(d\underline{x}) E^W([X_0(t)]^m \mid X(0) = \underline{x}) \\ &= E^{\xi^{(m)}}(\exp[bT_m(t)]\, f(\xi^1(t), \dots, \xi^m(t))), \end{aligned} \tag{2.57}$$

where the last line follows after substituting the representation formula in Lemma 2.1 [with an arbitrary initial condition $X(0)$] and doing a calculation similar to the one in the proof of Lemma 2.4. Passing to the limit $t \to \infty$ in (2.57), we get, with the help of (2.56), that

$$E^{\nu_\Theta}([X_0]^m) < \infty \quad \Longrightarrow \quad E^{\xi^{(m)}}(\exp[bT_m(\infty)]) < \infty. \tag{2.58}$$

Hence $\bar{b}_m \geq b_m$.

(d) We need to show that (2.54) is an equality. This is trivial when the right-hand side of (2.54) is infinite. Therefore, assume that $\nu_\Theta \neq \delta_{\underline{0}}$ and $E^{\nu_\Theta}([X_0]^m) < \infty$. By applying the mixing property of $\nu_\Theta$ [which was proved in part (b)] to (2.57), we get

$$E^{\nu_\Theta}([X_0]^m) = \Theta^m E^{\xi^{(m)}}(\exp[bT_m(\infty)]), \tag{2.59}$$

where we use (2.58), dominated convergence and the fact that $\xi^{(m)}$-a.s. all $m$ random walks move apart as $t \to \infty$ by transience. Moreover, by passing to the limit $t \to \infty$ in (2.18), we have

$$\lim_{t\to\infty} E([X_0(t)]^m \mid X(0) = \underline{\Theta}) = \Theta^m E^{\xi^{(m)}}(\exp[bT_m(\infty)]). \tag{2.60}$$

Combine (2.59)–(2.60) to get the claim.

2.5. *Properties of* $G^{(m)}$ *and* $K^{(m)}$. This section lists a number of elementary facts, many of which use basic random walk theory as explained in Spitzer [43].

Recall (1.36) and (1.39). We begin with the following statement:

LEMMA 2.6. *If* $a(\cdot,\cdot)$ *is symmetric and transient, then* $G^{(m)}(\cdot,\cdot)$ *is strongly transient for all* $m \geq 3$, *that is,*

$$\sup_{x,y\in I^{(m)}} \sum_{z\in I^{(m)}} G^{(m)}(x,z)G^{(m)}(z,y) < \infty, \qquad m \geq 3. \tag{2.61}$$

PROOF. Let

$$P_t^{(m)}(x,y) = P^{\eta^{(m)}}(\eta^{(m)}(t) = y \mid \eta^{(m)}(0) = x), \qquad x,y \in I^{(m)},\ t \geq 0, \tag{2.62}$$

be the transition probabilities of the differences random walk $\eta^{(m)} = (\eta^{(m)}(t))_{t\geq 0}$ defined in (1.33). We have

$$G^{(m)}(x,y) = m\int_0^\infty dt\, P_t^{(m)}(x,y), \qquad x,y \in I^{(m)}. \tag{2.63}$$

Compute

$$\begin{aligned} &\frac{1}{m^2}\sum_{z\in I^{(m)}} G^{(m)}(x,z)G^{(m)}(z,y) \\ &\quad= \int_0^\infty ds_1 \int_0^\infty ds_2 \sum_{z\in I^{(m)}} P_{s_1}^{(m)}(x,z)P_{s_2}^{(m)}(z,y) \\ &\quad= \int_0^\infty ds_1 \int_0^\infty ds_2\, P_{s_1+s_2}^{(m)}(x,y) \\ &\quad= \int_0^\infty dt\ tP_t^{(m)}(x,y) \\ &\quad\leq \int_0^\infty dt\ tP_t^{(m)}(0,0). \end{aligned} \tag{2.64}$$

Note that

$$P_t^{(m)}(0,0) = \sum_{i\in I}[P_t(0,i)]^m \tag{2.65}$$

with $P_t(i,j)$, $i,j \in I$, $t \geq 0$, the transition probabilities of a single random walk. Because the single random walk is symmetric and has exponential jump times (with mean 1), we have

$$\begin{aligned} &\text{(i)} \quad P_t(i,j) \leq P_t(0,0)\ \forall i,j \in I, t \geq 0, \\ &\text{(ii)} \quad t \mapsto P_t(0,0) \text{ is nonincreasing,} \end{aligned} \tag{2.66}$$

as is easily seen from the Fourier representation of $P_t(i,j)$; that is,

$$\begin{aligned} P_t(i,j) &= \sum_{n=0}^{\infty} e^{-t}\frac{t^n}{n!}\,(2\pi)^{-d}\int_{[-\pi,\pi)^d} d\lambda\,\cos(i-j,\lambda)[A(\lambda)]^n \\ &= (2\pi)^{-d}\int_{[-\pi,\pi)^d} d\lambda\,\cos(i-j,\lambda)e^{-t[1-A(\lambda)]} \end{aligned} \tag{2.67}$$

with $A(\lambda)=\sum_{i\in I}a(0,i)\cos(i,\lambda)\le 1$, $\lambda\in[-\pi,\pi)^d$, and $(\cdot,\cdot)$ the inner product on $\mathbb{R}^d$. Via (2.66)(i), (2.65) gives

$$P_t^{(m)}(0,0)\le [P_t(0,0)]^{m-2}\sum_{i\in I}P_t(0,i)P_t(i,0)=[P_t(0,0)]^{m-2}\,P_{2t}(0,0), \tag{2.68}$$

which via (2.65)(ii) yields

$$\begin{aligned} \int_0^\infty dt\;tP_t^{(m)}(0,0) &\le \tfrac{1}{2}\int_0^\infty dt\,[P_t(0,0)]^{m-2}\int_0^{2t} ds P_s(0,0) \\ &\le \tfrac{1}{2}\biggl(\int_0^\infty dt\,[P_t(0,0)]^{m-2}\biggr)\biggl(\int_0^\infty ds\,P_s(0,0)\biggr). \end{aligned} \tag{2.69}$$

The right-hand side is finite by the transience of the single random walk. □

We next look at $K^{(m)}$ defined in (1.39).

PROPOSITION 2.7. *Suppose that $a(\cdot,\cdot)$ is symmetric and transient. Then, for all $m\ge 2$, $K^{(m)}$ is a self-adjoint, positive and bounded operator on $\ell^2(S^{(m)})$.*

PROOF. The symmetry of $K^{(m)}$ follows from the symmetry of $G^{(m)}$. Since $K^{(m)}$ is defined everywhere on $\ell^2(S^{(m)})$, it therefore is self-adjoint. The Fourier representation of $G^{(m)}$ reads

$$\begin{aligned} G^{(m)}(x,y) &= \sum_{n=0}^{\infty}[a^{(m)}]^n(x,y) \\ &= |\widehat{I}^{(m)}|^{-1}\int_{\widehat{I}^{(m)}} d\widehat{\lambda}\;e^{i(x-y,\widehat{\lambda})}[1-A^{(m)}(\widehat{\lambda})]^{-1}, \qquad x,y\in I^{(m)}, \end{aligned} \tag{2.70}$$

with $\widehat{I}^{(m)}=([-\pi,\pi)^d)^{(1/2)m(m-1)}$, $(\cdot,\cdot)$ the inner product on $(\mathbb{R}^d)^{(1/2)m(m-1)}$, and $A^{(m)}(\widehat{\lambda})=\sum_{x\in I^{(m)}}a^{(m)}(0,x)\cos(x,\widehat{\lambda})\le 1$, $\widehat{\lambda}\in\widehat{I}^{(m)}$. It follows that

$$\begin{aligned} &\langle\mu,K^{(m)}\mu\rangle \\ &\quad= |\widehat{I}^{(m)}|^{-1}\int_{\widehat{I}^{(m)}} d\widehat{\lambda}\;[1-A^{(m)}(\widehat{\lambda})]^{-1}\biggl|\sum_{x\in I^{(m)}}e^{i(x,\widehat{\lambda})}\mu(x)\sqrt{\sharp^{(m)}(x)}\biggr|^2 \end{aligned} \tag{2.71}$$

with $\langle\cdot,\cdot\rangle$ the inner product on $\ell^2(S^{(m)})$. This proves the positivity of $K^{(m)}$. To prove the boundedness of $K^{(m)}$, we consider the relation

$$\|\sqrt{K^{(m)}}\|^2 = \sup_{\substack{\mu\in\ell^2(S^{(m)})\\ \langle\mu,\mu\rangle=1}} \langle\mu, K^{(m)}\mu\rangle \tag{2.72}$$

with $\|\cdot\|$ denoting the operator norm on $\ell^2(S^{(m)})$. Apply Cauchy–Schwarz twice, to obtain

$$\begin{aligned}
&\langle\mu, K^{(m)}\mu\rangle\\
&\quad= \sum_{x,y\in S^{(m)}} \mu(x)\sqrt{\sharp^{(m)}(x)}G^{(m)}(x,y)\sqrt{\sharp^{(m)}(y)}\mu(y)\\
&\quad\le \sum_{x\in S^{(m)}} \left[\sum_{y\in S^{(m)}} \sharp^{(m)}(x)G^{(m)}(x,y)\mu^2(y)\right]^{1/2}\\
&\qquad\qquad\times\left[\sum_{y\in S^{(m)}} \mu^2(x)G^{(m)}(x,y)\sharp^{(m)}(y)\right]^{1/2}\\
&\quad\le \left[\sum_{x,y\in S^{(m)}} \sharp^{(m)}(x)G^{(m)}(x,y)\mu^2(y)\right]^{1/2}\\
&\qquad\times\left[\sum_{x,y\in S^{(m)}} \mu^2(x)G^{(m)}(x,y)\sharp^{(m)}(y)\right]^{1/2}\\
&\quad\le \left[\sum_{x\in S^{(m)}} \mu^2(x)\right] \sup_{x\in S^{(m)}} \sum_{y\in S^{(m)}} G^{(m)}(x,y)\sharp^{(m)}(y),
\end{aligned} \tag{2.73}$$

where in the last line we use the symmetry of $G^{(m)}$. The last sum is equal to the average total collision local time (in pairs) of the $m$ walks when their differences start in $x$. Clearly, the supremum is taken at $x=0$, and equals $\sharp^{(m)}(0)G^{(2)}(0,0)$, because $G^{(2)}(0,0)$ is the average collision local time for each pair. Hence

$$\|\sqrt{K^{(m)}}\|^2 \le \sharp^{(m)}(0)G^{(2)}(0,0) < \infty. \tag{2.74}$$

But, by the self-adjointness and positivity of $K^{(m)}$, we have (see [39], Chapter 12)

$$\|\sqrt{K^{(m)}}\|^2 = \|[K^{(m)}]^n\|^{1/n} = \operatorname{spec}(K^{(m)}) \qquad \forall n\in\mathbb{N} \tag{2.75}$$

with $\operatorname{spec}(\cdot)$ denoting the spectral radius in $\ell^2(S^{(m)})$. $\square$

**3. Variational representations.** In Section 3.1 we identify $\bar{b}_m$ in terms of a variational formula. In Section 3.2 we prove that the $m$th moment diverges at $b = \bar{b}_m$. In Section 3.3 we calculate the exponential growth rate of the $m$th moment and prove that $\tilde{b}_m = \bar{b}_m$. In Section 3.4 we study the $m$-dependence of $b_m$. In Section 3.5 we collect the results and prove Theorem 1.6.

3.1. *Variational formula for* $\bar{b}_m$.

PROPOSITION 3.1. *Suppose that* $a(\cdot,\cdot)$ *is symmetric and transient. Then* $\bar{b}_m = m/\bar{\lambda}_m$ *with*

$$\bar{\lambda}_m = \sup_{\substack{\zeta\in\ell^1(S^{(m)})\\ \zeta\neq 0}} \frac{\langle \zeta, [K^{(m)}]^2\zeta\rangle}{\langle \zeta, K^{(m)}\zeta\rangle}. \tag{3.1}$$

PROOF. The proof comes in several steps. Throughout the proof we assume that $b\sharp^{(m)}(0) < m$.

*Step* 1. Recall the definition of $\bar{b}_m$ in (2.1) as well as the identity in (2.18). We begin by deriving a criterion for the property $E(\exp[bT^{(m)}(\infty)]) < \infty$ in terms of the discrete-time random walk

$$\eta^{(m),\odot} = (\eta^{(m),\odot}(i))_{i\in\mathbb{N}_0} \tag{3.2}$$

embedded in the continuous-time random walk $\eta^{(m)} = (\eta^{(m)}(t))_{t\geq 0}$ defined in (1.33). To that end we perform the expectation over the jump times of $\eta^{(m)}$, which are independent of $\eta^{(m),\odot}$. Indeed, let

$$(\sigma_i)_{i\in\mathbb{N}_0} \tag{3.3}$$

be the successive discrete times at which $\eta^{(m),\odot}$ visits $S^{(m)}$ (put $\sigma_0 = 0$), and let $M$ be its total number of visits to $S^{(m)}$ (which is random but finite a.s. by transience). Each visit to $S^{(m)}$ lasts a time $\tau$ that is exponentially distributed with mean $\frac{1}{m}$. Define

$$\square^{(m),b}(x) = E(\exp[b\sharp^{(m)}(x)\tau]) = \frac{m}{m - b\sharp^{(m)}(x)}, \qquad x\in S^{(m)}. \tag{3.4}$$

Then we have

$$E(\exp[bT^{(m)}(\infty)]) = E\left(\prod_{i=0}^{M} \square^{(m),b}(\eta^{(m),\odot}(\sigma_i))\right). \tag{3.5}$$

*Step* 2. In order to analyze the right-hand side of (3.5), we introduce the transition kernel of the Markov chain on $S^{(m)}$ obtained by observing $\eta^{(m),\odot}$ only when it visits $S^{(m)}$, which we denote by $P^{(m)}(\cdot,\cdot)$. Since $a(\cdot,\cdot)$ is symmetric, so is $P^{(m)}(\cdot,\cdot)$. By transience, this transition kernel is defective:

$$\sum_{y\in S^{(m)}} P^{(m)}(x,y)\begin{cases}\le 1, & \text{if } \sharp^{(m)}(x)=1,\\ =1, & \text{otherwise.}\end{cases} \tag{3.6}$$

(The first line says that escape from $S^{(m)}$ is possible only when all walks are disjoint except one pair. This is because only one walk moves at a time.) In terms of $P^{(m)}(\cdot,\cdot)$ we can write

$$\begin{aligned} &E\left(\prod_{i=0}^{M}\Box^{(m),b}(\eta^{(m),\odot}(\sigma_i))\right)\\ &\quad=\sum_{n=0}^{\infty}\sum_{x_0,\dots,x_n\in S^{(m)}}\delta_0(x_0)\left(\prod_{i=1}^{n}P^{(m)}(x_{i-1},x_i)\right)[1-P^{(m)}(x_n,S^{(m)})]\\ &\qquad\times\left(\prod_{i=0}^{n}\Box^{(m),b}(x_i)\right). \end{aligned} \tag{3.7}$$

Define the matrix

$$Q^{(m),b}(x,y)=\sqrt{\Box^{(m),b}(x)}P^{(m)}(x,y)\sqrt{\Box^{(m),b}(y)},\qquad x,y\in S^{(m)}. \tag{3.8}$$

With this notation we can write, combining (3.5) and (3.7)–(3.8),

$$E(\exp[bT^{(m)}(\infty)])=C_{m,b}\sum_{n=0}^{\infty}\langle\delta_0,[Q^{(m),b}]^nR_m\rangle \tag{3.9}$$

with

$$C_{m,b}=m/\sqrt{(m-b\sharp^{(m)}(0))(m-b)}$$

and

$$R_m(\cdot)=[1-P^{(m)}(\cdot,S^{(m)})]1_{\{\sharp^{(m)}(\cdot)=1\}}.$$

The front factor, which arises from the endpoint in the second sum in (3.7), is harmless.

*Step* 3. Note the following:

LEMMA 3.2. *$Q^{(m),b}(\cdot,\cdot)$ is an irreducible, aperiodic, nonnegative and symmetric matrix. As an operator acting on $\ell^2(S^{(m)})$ it is self-adjoint and bounded.*

PROOF. Because $\square^{(m),b}$ is bounded, $Q^{(m),b}(\cdot,\cdot)$ inherits these properties from $P^{(m)}(\cdot,\cdot)$. The irreducibility of $P^{(m)}(\cdot,\cdot)$ is inherited from the irreducibility of $a(\cdot,\cdot)$ assumed in (1.9). The aperiodicity of $P^{(m)}(\cdot,\cdot)$ follows from the fact that $P^{(m)}(x,x)>0$ for some $x\in S^{(m)}$ with $\sharp^{(m)}(x)=1$. □

Define

$$\bar\chi_m(b)=\lim_{n\to\infty}\frac{1}{n}\log[Q^{(m),b}]^n(x,y),\qquad x,y\in S^{(m)}. \tag{3.10}$$

Under the properties stated in Lemma 3.2, this limit exists, is in $\mathbb{R}$ and is the same for all $x,y\in S^{(m)}$ (see Vere-Jones [45]). Moreover,

$$\begin{aligned}&\sum_{n=0}^{\infty}[Q^{(m),b}]^n(0,0)<\infty\\ &\qquad\Longleftrightarrow\quad \sum_{n=0}^{\infty}[Q^{(m),b}]^n(x,y)<\infty\qquad \forall x,y\in S^{(m)}.\end{aligned} \tag{3.11}$$

This leads to

$$\begin{aligned}\bar\chi_m(b)>0 &\quad\Longrightarrow\quad \sum_{n=0}^{\infty}[Q^{(m),b}]^n(0,0)=\infty,\\ \bar\chi_m(b)<0 &\quad\Longrightarrow\quad \sum_{n=0}^{\infty}[Q^{(m),b}]^n(0,0)<\infty.\end{aligned} \tag{3.12}$$

From (2.1), (2.18), (3.5), (3.9) and (3.12) we see that $\bar b_m$ is the solution of the equation $\bar\chi_m(b)=0$. At the end of Section 3.2 the case $b=\bar b_m$ will be included in the top line of (3.12).

*Step* 4. Next, at $b=\bar b_m$ we have the following:

LEMMA 3.3.

$$\sup_{\substack{\nu\in\ell^1(S^{(m)})\\ \nu\neq 0}}\frac{\langle\nu,Q^{(m),\bar b_m}\nu\rangle}{\langle\nu,\nu\rangle}=1. \tag{3.13}$$

PROOF. Since $Q^{(m),\bar b_m}$ is a bounded operator on $\ell^2(S^{(m)})$, the function $\nu\mapsto\langle\nu,Q^{(m),\bar b_m}\nu\rangle$ is continuous on $\ell^2(S^{(m)})$. Since $\ell^1(S^{(m)})\subset\ell^2(S^{(m)})$ is dense, it suffices to prove that

$$\sup_{\substack{\nu\in\ell^2(S^{(m)})\\ \langle\nu,\nu\rangle=1}}\langle\nu,Q^{(m),\bar b_m}\nu\rangle=1. \tag{3.14}$$

[$\leq 1$]: First consider $\nu$ with finite support. Suppose that $\langle \nu, Q^{(m),\bar{b}_m}\nu\rangle \geq 1+\varepsilon$ for some $\varepsilon > 0$. Then, by the spectral theorem and Jensen's inequality,

$$\begin{aligned}\langle \nu, [Q^{(m),\bar{b}_m}]^{2n}\nu\rangle &= \int_{\mathbb{R}} \lambda^{2n}\, dE_{\nu,\nu}(\lambda) \geq \left(\int_{\mathbb{R}} \lambda\, dE_{\nu,\nu}(\lambda)\right)^{2n} \\ &= \langle \nu, Q^{(m),\bar{b}_m}\nu\rangle^{2n} \geq (1+\varepsilon)^{2n} \qquad \forall n \in \mathbb{N}\end{aligned} \tag{3.15}$$

with $E_{\nu,\nu}$ the spectral measure associated with $\nu$. Clearly this contradicts (3.10) with $\bar{\chi}_m(\bar{b}_m) = 0$, and so $\langle \nu, Q^{(m),\bar{b}_m}\nu\rangle \leq 1$ for $\nu$ with finite support. Since the $\nu$'s with finite support are dense in $\ell^2(S^{(m)})$, it follows that the supremum is $\leq 1$.

[$\geq 1$]: Suppose that the supremum is $\leq 1-\varepsilon$ for some $\varepsilon > 0$. Then the spectrum of $Q^{(m),\bar{b}_m}$ is contained in $(-\infty, 1-\varepsilon]$. Estimate

$$\begin{aligned}0 \leq \langle \delta_0, [Q^{(m),\bar{b}_m}]^{2n+1}\delta_0\rangle &= \int_{-\infty}^{1-\varepsilon} \lambda^{2n+1}\, dE_{\delta_0,\delta_0}(\lambda) \\ \leq \int_0^{1-\varepsilon} \lambda^{2n+1}\, dE_{\delta_0,\delta_0}(\lambda) &\leq (1-\varepsilon)^{2n+1} \qquad \forall n \in \mathbb{N}.\end{aligned} \tag{3.16}$$

But again this contradicts (3.10) with $\bar{\chi}_m(\bar{b}_m) = 0$. Hence the supremum is $\geq 1$. $\square$

*Step* 5. Putting $\mu = \sqrt{\square^{(m),\bar{b}_m}}\,\nu$, we may rewrite (3.13) as [recall (3.4)]

$$\begin{aligned}1 &= \sup_{\substack{\mu\in\ell^1(S^{(m)})\\ \mu\neq 0}} \frac{\langle \mu, P^{(m)}\mu\rangle}{\langle \mu, (\square^{(m),\bar{b}_m})^{-1}\mu\rangle} \\ &= \sup_{\substack{\mu\in\ell^1(S^{(m)})\\ \mu\neq 0}} \frac{\langle \mu,\mu\rangle - \langle \mu, (\mathbb{1}-P^{(m)})\mu\rangle}{\langle \mu,\mu\rangle - (\bar{b}_m/m)\langle \mu, \sharp^{(m)}\mu\rangle}.\end{aligned} \tag{3.17}$$

Therefore

$$\bar{b}_m = \frac{m}{\bar{\lambda}_m} \qquad \text{with } \bar{\lambda}_m = \sup_{\substack{\mu\in\ell^1(S^{(m)})\\ \mu\neq 0}} \frac{\langle \mu, \sharp^{(m)}\mu\rangle}{\langle \mu, (\mathbb{1}-P^{(m)})\mu\rangle}, \tag{3.18}$$

where the denominator is strictly positive because $P^{(m)}$ is irreducible. Let

$$\widehat{G}^{(m)} = \sum_{n=0}^{\infty} [P^{(m)}]^n = (\mathbb{1}-P^{(m)})^{-1}. \tag{3.19}$$

Because $P^{(m)}$ has spectral radius $< 1$ (due to the fact that $S^{(m)}$ is a uniformly transient set), we know that $\mathbb{1}-P^{(m)}$ is one-to-one on $\ell^1(S^{(m)})$ and that $\widehat{G}^{(m)}$ is a bounded operator on $\ell^1(S^{(m)})$. Therefore we can transform (3.18) via the change of variables $\mu = \widehat{G}^{(m)}\rho$:

$$\bar{\lambda}_m = \sup_{\substack{\rho\in\ell^1(S^{(m)})\\ \rho\neq 0}} \frac{\langle(\widehat{G}^{(m)}\rho)^2, \sharp^{(m)}\rangle}{\langle\rho, \widehat{G}^{(m)}\rho\rangle}. \tag{3.20}$$

Finally, putting $\rho = \sqrt{\sharp^{(m)}}\zeta$ and using that $\widehat{G}^{(m)}(x,y) = G^{(m)}(x,y)$ for all $x,y \in S^{(m)}$ by the definition of $P^{(m)}$, we get the formula in Proposition 3.1. □

3.2. *The $m$th moment at $b = \bar{b}_m$.* The case $b = \bar{b}_m$ can be included in the top line of (3.12) when 1 is the largest $\ell^1$-eigenvalue of $Q^{(m),\bar{b}_m}$. Therefore we next consider the eigenvalue equation

$$\nu Q^{(m),\bar{b}_m} = \nu, \qquad \nu \in \ell^1(S^{(m)}), \nu > 0. \tag{3.21}$$

LEMMA 3.4. *Suppose that $a(\cdot,\cdot)$ is symmetric and transient. If $\bar{b}_{m-1} > \bar{b}_m$, then (3.21) has a solution.*

PROOF. The idea is to use the notion of a quasi-stationary distribution.

*Step* 1. Consider the matrix

$$\begin{aligned} Q^{(m),\bar{b}_m,\otimes}(x,y) &= \frac{1}{N_{m,\bar{b}_m}} Q^{(m),\bar{b}_m}(x,y), \\ N_{m,\bar{b}_m} &= \sup_{x\in S^{(m)}} \sum_{y\in S^{(m)}} Q^{(m),\bar{b}_m}(x,y). \end{aligned} \tag{3.22}$$

This is an irreducible defective probability kernel on $S^{(m)}$. By introducing a cemetery state $\partial$, we can extend $Q^{(m),\bar{b}_m,\otimes}$ to a nondefective probability kernel on $S^{(m)} \cup \{\partial\}$. Let $(Z_n)_{n\in\mathbb{N}_0}$ denote the corresponding Markov chain starting in 0, and let

$$\nu_n = \mathcal{L}(Z_n \mid Z_n \neq \partial), \qquad n \in \mathbb{N}_0. \tag{3.23}$$

If we manage to show that ($\mathcal{P}$ denotes the set of probability measures)

$$\lim_{n\to\infty} \nu_n = \nu_\infty \qquad \text{in } \mathcal{P}(S^{(m)}), \tag{3.24}$$

then, because

$$\nu_{n+1} = \frac{\nu_n Q^{(m),\bar{b}_m,\otimes}}{\langle \nu_n Q^{(m),\bar{b}_m,\otimes}, 1\rangle}, \tag{3.25}$$

we get that

$$\nu_\infty Q^{(m),\bar{b}_m,\otimes} = \lambda_\infty \nu_\infty, \qquad \lambda_\infty > 0, \nu_\infty > 0, \tag{3.26}$$

with $\lambda_\infty = \langle \nu_\infty Q^{(m),\bar{b}_m,\otimes}, 1\rangle = 1/N_{m,\bar{b}_m}$ the probability of no defection to $\partial$ (per step) in the quasi-stationary distribution $\nu_\infty$. Hence $\nu_\infty$ solves (3.21).

*Step* 2. To prove (3.24), we use a criterion in Ferrari, Kesten and Martinez [23], Theorem 1, according to which it is enough to prove that there exist $\delta > 0$ and $D < \infty$ (depending on $m$) such that

$$P_0(\tau_0 > n \mid \tau_\partial > n) \le De^{-\delta n} \qquad \forall n \in \mathbb{N}_0 \tag{3.27}$$

with $P_0$ the law of $(Z_n)_{n\in\mathbb{N}_0}$ given $Z_0 = 0$, and $\tau_0, \tau_\partial$ the first hitting times of $0, \partial$ (time zero excluded).

For $K \subset \{1, \dots, m\}$ with $0 < |K| < m$, let

$$\begin{aligned} V_K = \{x \in S^{(m)} : &\text{for site } x \text{ there exist } i \in K, j \in K^c \\ &\text{such that walks } i \text{ and } j \text{ coincide}\} \end{aligned} \tag{3.28}$$

with $K^c = \{1, \dots, m\} \setminus K$ (recall from Section 1.6.1 that each $x \in I^{(m)}$ corresponds to a certain intersection order of the $m$ random walks). We will prove that there exist $\delta > 0$ and $D_K < \infty$ such that

$$P_0(\tau_{V_K} > n \mid \tau_\partial > n) \le D_K e^{-\delta n} \qquad \forall n \in \mathbb{N}_0, \tag{3.29}$$

where $\tau_{V_K}$ is the first hitting time of $V_K$ (time zero excluded). Since

$$\bigcap_{0<|K|<m} V_K = \{0\}, \tag{3.30}$$

we have $\{\tau_0 > n\} \subset \bigcup_{0<|K|<m} \{\tau_{V_K} > n\}$, and so (3.29) implies (3.27) with $D = \sum_K D_K$.

To prove (3.29), write

$$\begin{aligned} &P_0(\tau_{V_K} > n \mid \tau_\partial > n) = N/D \qquad \text{with} \\ &N = \sum_{y_1,\dots,y_n \in S^{(m)} \setminus V_K} Q^{(m),\bar{b}_m,\otimes}(0, y_1) Q^{(m),\bar{b}_m,\otimes}(y_1, y_2) \times \cdots \\ &\qquad\qquad \times Q^{(m),\bar{b}_m,\otimes}(y_{n-1}, y_n), \\ &D = \sum_{y_1,\dots,y_n \in S^{(m)}} Q^{(m),\bar{b}_m,\otimes}(0, y_1) Q^{(m),\bar{b}_m,\otimes}(y_1, y_2) \times \cdots \\ &\qquad\qquad \times Q^{(m),\bar{b}_m,\otimes}(y_{n-1}, y_n) \end{aligned} \tag{3.31}$$

We may drop the $\otimes$ because the normalization factor in (3.22) cancels out. After that the denominator in (3.31) equals

$$\langle \delta_0, [Q^{(m),\bar{b}_m}]^n 1\rangle \geq \langle \delta_0, [Q^{(m),\bar{b}_m}]^n \delta_0\rangle = [\bar{\chi}_m(\bar{b}_m) + o(1)]^n = \exp[o(n)] \tag{3.32}$$

because $\bar{\chi}_m(\bar{b}_m) = 0$. It therefore suffices to prove that the numerator in (3.31) satisfies the exponential bound in (3.29).

Now, on $S^{(m)} \setminus V_K$ we have

$$T^{(m)}(\sigma_n) = T^{(m),K}(\sigma_n) + T^{(m),K^c}(\sigma_n) \qquad \forall n \in \mathbb{N} \tag{3.33}$$

with $\sigma_n$ the time of the $n$th visit to $S^{(m)}$ [recall (3.3)] and $T^{(m),K}(\sigma_n)$ the total collision local time (in pairs) up to time $\sigma_n$ of the walks indexed by $K$, and similarly for $K^c$. Therefore, retracing the calculations in Steps 1 and 2 of the proof of Proposition 3.1, we find that

$$\text{numerator (3.31)} \leq E_0(\exp[\bar{b}_m\{T^{(m),K}(\sigma_n) + T^{(m),K^c}(\sigma_n)\}]1_{\{\sigma_n<\infty\}}). \tag{3.34}$$

The inequality comes from using (3.33) and afterward dropping the restriction to $S^{(m)} \setminus V_K$. Next, apply Hölder's inequality to get, for $\varepsilon > 0$,

$$\begin{aligned} &\text{numerator (3.31)} \\ &\qquad \leq E_0(\exp[\bar{b}_m\{T^{(m),K}(\infty) + T^{(m),K^c}(\infty)\}]1_{\{\sigma_n<\infty\}}) \\ &\qquad \leq E_0(\exp[(1+\varepsilon)\bar{b}_m\{T^{(m),K}(\infty) + T^{(m),K^c}(\infty)\}])^{1/(1+\varepsilon)} \\ &\qquad\quad \times P_0(\sigma_n < \infty)^{\varepsilon/(1+\varepsilon)}. \end{aligned} \tag{3.35}$$

The expectation in the right-hand side factors because $T^{(m),K}(\infty)$ and $T^{(m),K^c}(\infty)$ are independent, and each factor is finite when $\varepsilon$ is picked so small that $(1+\varepsilon)\bar{b}_m < \widehat{b}_{m-1}$, because $|K|, |K^c| \leq m-1$. On the other hand, the probability in the right-hand side tends to zero exponentially fast with $n$ because $S^{(m)}$ is a uniformly transient set. $\square$

By (2.1), we have $E^{\nu_\Theta}([X_0]^m) < \infty$ for $b < b_m$ and $E^{\nu_\Theta}([X_0]^m) = \infty$ for $b > b_m$. With the help of Lemma 3.4 we can now include $b = b_m$.

LEMMA 3.5. *$E^{\nu_\Theta}([X_0]^m) = \infty$ at $b = b_m$.*

PROOF. By Proposition 2.5(b), we have $b_m = \bar{b}_m$. Suppose first that $b_{m-1} > b_m$. Then, by Lemma 3.4,

$$\nu Q^{(m),b_m} = \nu \qquad \nu \in \ell^1(S^{(m)}), \text{for some } \nu > 0, \tag{3.36}$$

and hence

$$\sum_{n=0}^{\infty} \sum_{x\in S^{(m)}} \nu(x)[Q^{(m),b_m}]^n(x,y) = \sum_{n=0}^{\infty} \nu(y) = \infty \qquad \forall y \in S^{(m)}. \tag{3.37}$$

By the irreducibility of $Q^{(m),b_m}$, this implies that

$$\sum_{n=0}^{\infty}[Q^{(m),b_m}]^n(0,0)=\infty, \tag{3.38}$$

which shows, via (3.5) and (3.9), that

$$E^{\xi^{(m)}}(\exp[b_m T^{(m)}(\infty)])=\infty. \tag{3.39}$$

Now use Proposition 2.5(d) to obtain from (3.39) that $E^{\nu_\Theta}([X_0]^m)=\infty$ at $b=b_m$.

Suppose next that $b_{m-1}=b_m$, but $b_{m'-1}>b_m$ for some $m'<m$. We may assume that $m'$ is the largest such index. Then $b_{m'-1}>b_{m'}$, and so the above argument tells us that $E^{\nu_\Theta}([X_0]^{m'})=\infty$ at $b=b_{m'}$. But because $b_{m'}=b_m$ and $m'<m$, we again get $E^{\nu_\Theta}([X_0]^m)=\infty$ at $b=b_m$.

Finally, if $b_{m'}=b_m$ for all $m'<m$, then $b_2=b_m$. However, at $b=b_2$ we have $E^{\nu_\Theta}([X_0]^2)=\infty$, as is easily seen from Lemma 2.4, Proposition 2.5(d) and (3.9), because $S^{(2)}=\{0\}$. Therefore once again $E^{\nu_\Theta}([X_0]^m)=\infty$ at $b=b_m$. $\square$

3.3. *Growth rate of the mth moment and* $\widetilde{b}_m=\bar{b}_m$. In this section we show that

$$\lim_{t\to\infty}\frac{1}{t}\log E^{\xi^{(m)}}(\exp[bT^{(m)}(t)])=\chi_m(b) \tag{3.40}$$

exists and can be expressed in terms of a variational problem. We will analyze this variational problem and show that

$$\widetilde{b}_m=\sup\{b>0:\chi_m(b)=0\} \tag{3.41}$$

[recall (2.1) and (2.18)] coincides with $\bar{b}_m$. Together with Proposition 2.5(c) this will show that all three critical values in (2.1) coincide.

In order to pose the problem in a form suitable for a large deviation analysis, we recall the definition of the differences random walk $(\eta^{(m)}(t))_{t\geq 0}$ in (1.33) and the collision function $\sharp^{(m)}$ in (1.37). Using (2.19), we have the identity

$$T^{(m)}(t)=\int_0^t \sharp^{(m)}(\eta^{(m)}(s))\,ds. \tag{3.42}$$

Next, we define the empirical measure

$$L_t^{(m)}=\frac{1}{t}\int_0^t \delta_{\eta^{(m)}(s)}\,ds. \tag{3.43}$$

LEMMA 3.6. *Suppose that $a(\cdot,\cdot)$ is symmetric. Then $(L_t^{(m)})_{t\ge 0}$ satisfies the weak large deviation principle on $\mathcal{P}(I^{(m)})$ with rate function*

$$J^{(m)}(\nu) = \langle \nu^{1/2}, m(\mathbb{1} - a^{(m)})\nu^{1/2}\rangle, \qquad \nu \in \mathcal{P}(I^{(m)}), \tag{3.44}$$

*where $a^{(m)}(\cdot,\cdot)$ is the transition kernel defined in (1.34). $J^{(m)}$ is bounded on $\ell^2(I^{(m)})$ and when restricted to $\ell^1(I^{(m)})$ is continuous in the $\ell^1$-topology.*

PROOF. See Deuschel and Stroock [21], Section 4.2. The rate function is given by (3.44) because $m(a^{(m)} - \mathbb{1})$ is the generator of the Markov process $(\eta^{(m)}(t))_{t\ge 0}$ and $a^{(m)}(\cdot,\cdot)$ is symmetric [recall (1.32)]. The latter is crucial for having the explicit formula in (3.44). The boundedness of $J^{(m)}$ is obvious. The continuity of $J^{(m)}$ follows from the fact that $m(\mathbb{1} - a^{(m)})$ is a bounded operator from $\ell^2(I^{(m)})$ into $\ell^1(I^{(m)})$, in combination with the fact that if $\nu_n \to \nu$ as $n \to \infty$ in $\ell^1$-norm, then $\nu_n^{1/2} \to \nu^{1/2}$ in $\ell^2$-norm. To see the latter, write

$$\|\nu_n^{1/2} - \nu^{1/2}\|_2^2 = \sum_x \bigl(\sqrt{\nu_n(x)} - \sqrt{\nu(x)}\bigr)^2 = 2 - 2\sum_x \sqrt{\nu_n(x)\nu(x)}. \tag{3.45}$$

If $\nu_n \to \nu$ in $\ell^1$-norm, then Fatou's lemma gives $\liminf_{n\to\infty} \sum_x \sqrt{\nu_n(x)\nu(x)} \ge \sum_x \nu(x) = 1$. Hence $\lim_{n\to\infty} \|\nu_n^{1/2} - \nu^{1/2}\|_2^2 = 0$. $\square$

Lemma 3.6 leads us to the following identification:

LEMMA 3.7.

$$\chi_m(b) = \sup_{\nu\in\mathcal{P}(I^{(m)})} \{b\langle \nu, \sharp^{(m)}\rangle - J^{(m)}(\nu)\}. \tag{3.46}$$

PROOF. For ease of notation we drop the superscript $(m)$. We cannot apply Varadhan's lemma directly to (3.40)–(3.43), since we only have a *weak* large deviation principle. This problem can be handled via a standard compactification argument, as follows.

Let $(\eta_N^+(t))_{t\ge 0}, (\eta_N^-(t))_{t\ge 0}$ be the differences random walks obtained by wrapping $(\eta(t))_{t\ge 0}$ around the torus $\Lambda_N = ([-N,N)^d \cap I)^{(1/2)m(m-1)}$ (recall that $I = \mathbb{Z}^d$), respectively, by killing it when it hits $\partial\Lambda_N$, the boundary of $\Lambda_N$. Let $T_N^+(t), T_N^-(t)$ be the quantities corresponding to $T(t)$ in (3.42) for these two processes. Then, by the large deviation principle in Lemma 3.6 restricted to $\Lambda_N$, we have for every $N$:

$$\begin{aligned} &\lim_{t\to\infty} \frac{1}{t}\log E(\exp[bT_N^+(t)]) = S_N^+,\\ &\lim_{t\to\infty} \frac{1}{t}\log E(\exp[bT_N^-(t)]) = S_N^-, \end{aligned} \tag{3.47}$$

with

$$S_N^+ = \sup_{\nu_N \in \mathcal{P}(\Lambda_N)} \{b\langle \nu_N, \sharp\rangle - J_N(\nu_N)\},$$

$$S_N^- = \sup_{\substack{\nu_N \in \mathcal{P}(\Lambda_N) \\ \nu_N(\partial\Lambda_N)=0}} \{b\langle \nu_N, \sharp\rangle - J_N(\nu_N)\}. \tag{3.48}$$

Here, $J_N$ is the analogue of $J$ in (3.44) restricted to $\Lambda_N$, that is,

$$J_N(\nu_N) = \langle \nu_N^{1/2}, m(\mathbb{1} - a_N)\nu_N^{1/2}\rangle, \qquad \nu_N \in \mathcal{P}(\Lambda_N), \tag{3.49}$$

with $a_N(\cdot,\cdot)$ the periodized transition kernel

$$a_N(x,y) = \sum_{z \in (2N)I^{(m)}} a(x, y+z), \qquad x, y \in \Lambda_N. \tag{3.50}$$

Now, obviously

$$T_N^+(t) \geq T(t) \geq T_N^-(t), \tag{3.51}$$

because the intersection local time increases by wrapping and decreases by killing. Consequently,

$$S_N^+ \geq S \geq S_N^- \tag{3.52}$$

with $S$ the right-hand side of (3.46). Hence, to prove (3.40) and (3.46) it suffices to show that

$$\liminf_{N\to\infty} S_N^- \geq S, \qquad \limsup_{N\to\infty} S_N^+ \leq S. \tag{3.53}$$

*Lower bound.* It suffices to show that for every $\nu \in \mathcal{P}(I^{(m)})$ there exists a sequence $(\nu_N)$, with $\nu_N \in \mathcal{P}(\Lambda_N)$ and $\nu_N(\partial\Lambda_N) = 0$ for each $N$, such that

$$\text{w-lim}_{N\to\infty} \nu_N = \nu, \qquad \lim_{N\to\infty} J_N(\nu_N) = J(\nu). \tag{3.54}$$

Indeed, with the help of Fatou's lemma this gives

$$\liminf_{N\to\infty} S_N^- \geq \liminf_{N\to\infty} \{b\langle \nu_N, \sharp\rangle - J_N(\nu_N)\} \geq b\langle \nu, \sharp\rangle - J(\nu), \tag{3.55}$$

and so we get the first half of (3.53) after taking the supremum over $\nu$ afterward.

For the given $\nu$, the sequence $(\nu_N)$ is chosen as follows. Put $\nu_N(x) = \nu(x)$ for all $x \in \Lambda_N \setminus (\partial\Lambda_N \cup \{0\})$ and $\nu_N(0) = \nu(0) + \sum_{x \notin \Lambda_N \cup \partial\Lambda_N} \nu(x)$. Then the first half of (3.54) is obvious. For the second half, since $J$ is continuous in the $\ell_1$-topology it suffices to show that $\lim_{N\to\infty}[J_N(\nu_N) - J(\nu_N)] = 0$. To that end, we estimate

$$
\begin{aligned}
0 &\le J(\nu_N) - J_N(\nu_N) \\
&= m\langle \nu_N^{1/2}, (a_N - a)\nu_N^{1/2}\rangle \\
&\le m \sum_{x,y\in\Lambda_N} \nu_N(x)(a_N - a)(x,y) \\
&\le m \sup_{x\in\Lambda_N} \sum_{y\in\Lambda_N} (a_N - a)(x,y) \\
&\le m \sup_{x\in\Lambda_N} \sum_{\substack{y\in I^{(m)} \\ \|y-x\|_\infty \ge N}} a(x,y) = \delta_N,
\end{aligned}
\tag{3.56}
$$

where in the third line we use the Cauchy–Schwarz inequality and the symmetry of $a_N$ and $a$, and in the fifth line we use (3.50) and the shift-invariance of $a(\cdot,\cdot)$ (note that on the sublattice with spacing $2N$ containing $y \in \Lambda_N$, the site closest to $x$ is captured by the supremum over $x \in \Lambda_N$ and the sum over $y \in I^{(m)}$ with $\|y - x\|_\infty < N$). Obviously, $\lim_{N\to\infty} \delta_N = 0$, which completes the proof of the first half of (3.53).

*Upper bound.* Estimate, with the help of (3.56),

$$
\begin{aligned}
S_N^+ &= \sup_{\nu_N\in\mathcal{P}(\Lambda_N)} \{b\langle \nu_N, \sharp\rangle - J_N(\nu_N)\} \\
&\le \sup_{\nu_N\in\mathcal{P}(\Lambda_N)} \{b\langle \nu_N, \sharp\rangle - J(\nu_N)\} + \delta_N \\
&\le S + \delta_N.
\end{aligned}
\tag{3.57}
$$

Let $N \to \infty$ to obtain the second half of (3.53). □

It follows from (3.46) that $b \mapsto \chi_m(b)$ is nondecreasing and convex on $[0,\infty)$. Since it is finite, it is also continuous on $[0,\infty)$. Obviously, $\chi_m(0) = 0$. The critical value $\widetilde{b}_m$ is given by (3.41). It further follows from (3.46) that

$$
b\,\sharp^{(m)}(0) - J^{(m)}(\delta_0) \le \chi_m(b) \le b\,\sharp^{(m)}(0), \tag{3.58}
$$

where $J^{(m)}(\delta_0) = m(1 - a^{(m)}(0,0)) = m$ [use (3.44) and (1.34), and recall that $a(0,0) = 0$ as assumed below (1.9)].

PROPOSITION 3.8. *Suppose that $a(\cdot,\cdot)$ is symmetric and transient. Then* $\widetilde{b}_m = \bar{b}_m$.

PROOF. Changing variables in (3.46) by putting $\nu = \mu^2$, we get

$$\chi_m(b) = \sup_{\substack{\mu \in \ell^2(I^{(m)}) \\ \langle \mu, \mu \rangle = 1}} F^{(m)}(\mu) \tag{3.59}$$

with

$$F^{(m)}(\mu) = b\langle \mu^2, \sharp^{(m)} \rangle - m\langle \mu, (\mathbb{1} - a^{(m)})\mu \rangle. \tag{3.60}$$

Define

$$\widetilde{\lambda}_m = \sup_{\substack{\mu \in \ell^2(I^{(m)}) \\ \mu \neq 0}} \frac{\langle \mu^2, \sharp^{(m)} \rangle}{\langle \mu, (\mathbb{1} - a^{(m)})\mu \rangle}, \tag{3.61}$$

where the denominator is strictly positive because $a^{(m)}(\cdot,\cdot)$ is irreducible. It follows from (3.59)–(3.60) that

$$\begin{aligned} b\widetilde{\lambda}_m > m \quad &\Longrightarrow \quad \chi_m(b) > 0, \\ b\widetilde{\lambda}_m \le m \quad &\Longrightarrow \quad \chi_m(b) = 0. \end{aligned} \tag{3.62}$$

Hence $\widetilde{b}_m = m/\widetilde{\lambda}_m$. Change variables in (3.61) by putting $\mu = G^{(m)}\rho$. Then, because $G^{(m)} = (\mathbb{1} - a^{(m)})^{-1}$, we get

$$\widetilde{\lambda}_m \ge \sup_{\substack{\rho \in \ell^1(S^{(m)}) \\ \rho \neq 0}} \frac{\langle (G^{(m)}\rho)^2, \sharp^{(m)} \rangle}{\langle \rho, G^{(m)}\rho \rangle}. \tag{3.63}$$

Here the inequality arises because we restrict the support of $\rho$ to $S^{(m)}: G^{(m)}\rho$ does not run through all of $\ell^2(S^{(m)})$. Note here that, because $G^{(m)}$ is strongly transient by Lemma 2.6, $\rho \in \ell^1(S^{(m)})$ implies that $\mu = G^{(m)}\rho \in \ell^2(I^{(m)})$.

Putting $\rho = \sqrt{\sharp^{(m)}}\zeta$ with $\zeta \in \ell^1(S^{(m)})$, we obtain

$$\widetilde{\lambda}_m \ge \sup_{\substack{\zeta \in \ell^1(S^{(m)}) \\ \zeta \neq 0}} \frac{\langle \zeta, [K^{(m)}]^2 \zeta \rangle}{\langle \zeta, K^{(m)}\zeta \rangle}. \tag{3.64}$$

Comparing with (3.1), we have thus found that $\widetilde{\lambda}_m \ge \bar{\lambda}_m$. Consequently, $\widetilde{b}_m \le \bar{b}_m$. However, it is obvious from (2.1) that $\widetilde{b}_m \ge \bar{b}_m$. Hence we get $\widetilde{b}_m = \bar{b}_m$ and $\widetilde{\lambda}_m = \bar{\lambda}_m$. □

3.4. *Analysis of* $m \mapsto b_m$. The remaining steps in this section are the following three propositions.

PROPOSITION 3.9. *Suppose that* $a(\cdot,\cdot)$ *is symmetric and transient. Then* $b_m = \bar{b}_m = \widetilde{b}_m = \frac{m}{\lambda_m}$ *with* $\lambda_m$ *the spectral radius of* $K^{(m)}$ *in* $\ell^2(S^{(m)})$.

PROOF. It follows from Propositions 2.5(b), 3.1 and 3.8 that $b_m = \bar{b}_m = \widetilde{b}_m = \frac{m}{\lambda_m}$ with $\lambda_m$ given by the right-hand side of (3.1). We saw in the proof of Lemma 3.3 that the supremum may be taken over $\nu \in \ell^2(S^{(m)})$. Putting $\mu = \sqrt{K^{(m)}}\nu$ in (3.1), we see that $\lambda_m = \|\sqrt{K^{(m)}}\|^2$. Now recall (2.75). $\square$

PROPOSITION 3.10. *Suppose that $a(\cdot,\cdot)$ is symmetric and transient. Then the spectral radius of $K^{(m)}$ is an eigenvalue if and only if $b_{m-1} > b_m$.*

PROOF. The fact that $b_{m-1} > b_m$ implies that the spectral radius of $K^{(m)}$ is an eigenvalue follows from Lemma 3.4 via the change of variables in Step 5 in the proof of Proposition 3.1. Indeed, Lemma 3.4 shows that the supremum in (3.1) is then actually attained.

To get the reverse, we use Lemma 3.4 and Proposition 3.9. The idea is to imbed the variational problem with label $m-1$ into the one with label $m$. To that end, consider $m$ random walks but ignore the collision local time the $m$th random walk has with the others. By repeating the large deviation analysis in Section 3.3, we get a formula for $\widetilde{\chi}_{m-1}(b)$ of the same form as in (3.46) but with the first term $\langle \nu, \sharp^{(m)} \rangle$ replaced by

$$\langle \nu, \sharp^{(m),(m-1)} \rangle, \tag{3.65}$$

where [cf. with (1.37)]

$$\begin{aligned} \sharp^{(m),(m-1)}(z) &= \sum_{1 \le p < q \le m-1} 1_{\{z_p - z_q = 0\}}, \\ & z = (z_p - z_q)_{1 \le p < q \le m},\ z_p, z_q \in I. \end{aligned} \tag{3.66}$$

In this way we can represent $\chi_m(b)$ and $\chi_{m-1}(b)$ via a variational problem on the same space:

$$\begin{aligned} \chi_m(b) &= \sup_{\nu \in \mathcal{P}(I^{(m)})} \{ b\langle \nu, \sharp^{(m)} \rangle - J^{(m)}(\nu) \}, \\ \chi_{m-1}(b) &= \sup_{\nu \in \mathcal{P}(I^{(m)})} \{ b\langle \nu, \sharp^{(m),(m-1)} \rangle - J^{(m)}(\nu) \}. \end{aligned} \tag{3.67}$$

Now consider these two variational problems at $b = b_m$. Suppose that $K^{(m)}$ has a maximal eigenvalue. Then this eigenvalue is unique, and so is its corresponding eigenvector. Consequently, the first variational problem has a unique maximizer, say $\bar{\nu} \in \mathcal{P}(I^{(m)})$, which is strictly positive. For $\varepsilon \in (0, \bar{\nu}(0))$, put

$$\mathcal{U}_\varepsilon = \{ \nu \in \mathcal{P}(I^{(m)}) : |\nu(0) - \bar{\nu}(0)| < \varepsilon \}, \tag{3.68}$$

and let

$$A_\varepsilon = \sup_{\nu \in \mathcal{P}(I^{(m)}) \setminus \mathcal{U}_\varepsilon} \{b_m \langle \nu, \sharp^{(m),(m-1)} \rangle - J^{(m)}(\nu)\},$$

(3.69)

$$B_\varepsilon = \sup_{\nu \in \mathcal{U}_\varepsilon} \{b_m \langle \nu, \sharp^{(m),(m-1)} \rangle - J^{(m)}(\nu)\}.$$

Since $\bar{\nu}$ is the unique maximizer, we have

$$0 = \chi_m(b_m) > \sup_{\nu \in \mathcal{P}(I^{(m)}) \setminus \mathcal{U}_\varepsilon} \{b_m \langle \nu, \sharp^{(m)} \rangle - J^{(m)}(\nu)\} \geq A_\varepsilon. \tag{3.70}$$

Since $\sharp^{(m)}(0) > \sharp^{(m),(m-1)}(0)$, we have

$$0 = \chi_m(b_m) \geq \sup_{\nu \in \mathcal{U}_\varepsilon} \{b_m \langle \nu, \sharp^{(m)} \rangle - J^{(m)}(\nu)\} > B_\varepsilon. \tag{3.71}$$

Combining (3.70)–(3.71) with the observation that $\chi_{m-1}(b_m) = A_\varepsilon \vee B_\varepsilon$, we find that $\chi_{m-1}(b_m) < 0$ and hence that $b_{m-1} > b_m$. □

PROPOSITION 3.11. *Suppose that $a(\cdot,\cdot)$ is transient. Then:*

(a) $b_2 = 2/G^{(2)}(0,0)$.

*Suppose that $a(\cdot,\cdot)$ is symmetric and transient. Then:*

(b) $b_2 \geq b_3 \geq b_4 \geq \cdots > 0$.

(c) $b_2 \leq (m-1)b_m \leq 2/G^{(m)}(0,0) < 2$.

(d) $\lim_{m\to\infty}(m-1)b_m = c = \sup_{m\geq 2}(m-1)b_m$.

PROOF. (a) The formula for $b_2$ is obvious because $S^{(2)} = \{0\}$ and $\sharp^{(2)}(0) = 1$, giving $\lambda_2 = G^{(2)}(0,0)$.

(b) The fact that $m \mapsto b_m$ is nonincreasing is trivial from (2.1) and Lemma 2.4. The fact that $b_m > 0$ for all $m$ follows from (c).

(c) We prove that

$$\sharp^{(m)}(0) G^{(m)}(0,0) \leq \lambda_m \leq \sharp^{(m)}(0) G^{(2)}(0,0), \tag{3.72}$$

which implies the claim. The lower bound is obtained by picking $\mu = \delta_0$ in (2.72) and using (2.75). The upper bound follows from (2.74)–(2.75).

(d) The proof is an adaptation of the argument in Carmona and Molchanov [7], Chapter III, Section 1.3. The key is the following inequality:

LEMMA 3.12. *If*

$$m = \sum_{i=1}^{r} n_i m_i \qquad \text{with } n_i, m_i \in \mathbb{N} \text{ and } m_i \geq 2 \text{ for } i = 1, \ldots, r, \tag{3.73}$$

*then*

$$\chi_m\left(\frac{b}{m-1}\right) \leq \sum_{i=1}^{r} n_i \chi_{m_i}\left(\frac{b}{m_i - 1}\right) \qquad \forall b > 0. \tag{3.74}$$

Proof. Let $\mathcal{P}$ denote the collection of all partitions of $\{1,\dots,m\}$ into $n_i$ groups of $m_i$ integers for $i=1,\dots,r$. For $P\in\mathcal{P}$, write

$$P=(P_{ij})_{i=1,\dots,r,j=1,\dots,n_i} \tag{3.75}$$

to label these groups, so that

$$\{1,\dots,m\}=\bigcup_{i=1}^{r}\bigcup_{j=1}^{n_i}P_{ij} \tag{3.76}$$

with

$$|P_{ij}|=m_i,\qquad i=1,\dots,r,\ j=1,\dots,n_i. \tag{3.77}$$

The cardinality of $\mathcal{P}$ is

$$|\mathcal{P}|=N=\frac{m!}{\prod_{i=1}^{r}n_i!\,m_i^{n_i}}. \tag{3.78}$$

Moreover,

$$\left[\sum_{P\in\mathcal{P}}\sum_{j=1}^{n_i}\sum_{\substack{k,l\in P_{ij}\\1\le k<l\le m}}1\right]=N_i,\qquad i=1,\dots,r, \tag{3.79}$$

with

$$N_i=n_i\frac{m_i(m_i-1)}{m(m-1)}N. \tag{3.80}$$

Now, let

$$\delta_i=\frac{m-1}{m_i-1}\,\frac{1}{N}. \tag{3.81}$$

Then, by (3.73),

$$\sum_{i=1}^{r}\delta_iN_i=1. \tag{3.82}$$

Using (3.79) and (3.82), we may write

$$\begin{aligned}T^{(m)}(t)&=\sum_{1\le k<l\le m}\int_0^t 1_{\{\xi_k(s)=\xi_l(s)\}}\,ds\\&=\sum_{P\in\mathcal{P}}\sum_{i=1}^{r}\delta_i\sum_{j=1}^{n_i}\sum_{\substack{k,l\in P_{ij}\\1\le k<l\le m}}\int_0^t 1_{\{\xi_k(s)=\xi_l(s)\}}\,ds.\end{aligned} \tag{3.83}$$

With the help of (3.83), we may estimate

$$
\begin{aligned}
&E^{\xi^{(m)}}(\exp[bT^{(m)}(t)]) \\
&\quad = E^{\xi^{(m)}}\Bigg(\prod_{P\in\mathcal{P}}\prod_{i=1}^{r}\exp\Bigg[b\delta_i\sum_{j=1}^{n_i}\sum_{\substack{k,l\in P_{ij}\\1\le k<l\le m}}\int_0^t 1_{\{\xi_k(s)=\xi_l(s)\}}\,ds\Bigg]\Bigg) \\
&\quad \le \prod_{P\in\mathcal{P}} E^{\xi^{(m)}}\Bigg(\prod_{i=1}^{r}\exp\Bigg[Nb\delta_i\sum_{j=1}^{n_i}\sum_{\substack{k,l\in P_{ij}\\1\le k<l\le m}}\int_0^t 1_{\{\xi_k(s)=\xi_l(s)\}}\,ds\Bigg]\Bigg)^{1/N} \qquad (3.84)\\
&\quad = \prod_{P\in\mathcal{P}}\prod_{i=1}^{r} E^{\xi^{(m)}}\Bigg(\exp\Bigg[Nb\delta_i\sum_{j=1}^{n_i}\sum_{\substack{k,l\in P_{ij}\\1\le k<l\le m}}\int_0^t 1_{\{\xi_k(s)=\xi_l(s)\}}\,ds\Bigg]\Bigg)^{1/N} \\
&\quad = \prod_{i=1}^{r} E^{\xi^{(m)}}\Bigg(\exp\Bigg[Nb\delta_i\sum_{j=1}^{n_i}\sum_{\substack{k,l\in P^*_{ij}\\1\le k<l\le m}}\int_0^t 1_{\{\xi_k(s)=\xi_l(s)\}}\,ds\Bigg]\Bigg),
\end{aligned}
$$

where $P^*$ is any representative partition. Here, the third line uses Hölder's inequality, the fourth line uses that the factors labelled by $i$ are independent, and the fifth line uses that the expectations for given $P$ do not depend on the choice of $P$. Taking logarithms, dividing by $t$ and letting $t\to\infty$ on both sides of (3.84), we obtain

$$\chi_m(b) \le \sum_{i=1}^{r} n_i\chi_{m_i}(Nb\delta_i) \qquad (3.85)$$

[recall (3.40) and (3.77)]. Inserting (3.81) and replacing $b$ by $b/(m-1)$, we arrive at (3.74). □

We now show why Lemma 3.12 implies the claim in Proposition 3.11(d). An immediate consequence of Lemma 3.12 and the definition of $b_m$,

$$b_m = \sup\{b>0 : \chi_m(b)=0\} \qquad (3.86)$$

[recall (3.41) and use that $b_m=\bar{b}_m$], is that, subject to (3.73),

$$(m-1)b_m \ge \min_{i=1,\dots,r}(m_i-1)b_{m_i}. \qquad (3.87)$$

Let $c=\limsup_{m\to\infty}(m-1)b_m$. Pick $\varepsilon>0$ arbitrary, and pick $m_1=m_1(\varepsilon)$ such that

$$(m_1-1)b_{m_1} \ge c-\varepsilon. \qquad (3.88)$$

For $m > m_1$, let $n_1 = \lceil m/m_1 \rceil$. Since $m \le n_1 m_1 < m + m_1$, we may estimate

$$(3.89)\quad \begin{aligned}(m-1)b_m &\ge (m-1)b_{n_1 m_1} = \frac{m-1}{n_1 m_1 - 1}(n_1 m_1 - 1)b_{n_1 m_1}\\ &\ge \frac{m-1}{n_1 m_1 - 1}(m_1 - 1)b_{m_1} > \frac{m-1}{m + m_1 - 1}(c - \varepsilon).\end{aligned}$$

Here, the first inequality uses that $m \mapsto b_m$ is nonincreasing, while the second inequality uses (3.87) for $r = 1$. Let $m \to \infty$ to obtain that $\liminf_{m\to\infty}(m-1)b_m \ge c - \varepsilon$. Let $\varepsilon \downarrow 0$ to get $\lim_{m\to\infty}(m-1)b_m = c$. Note that $c = \sup_{m\ge 2}(m-1)b_m$. $\square$

3.5. *Proof of Theorem* 1.6. (a) Combine (2.1), Proposition 2.5(d) and Lemma 3.5.

(b) Combine (2.18), (3.40), (3.41), Propositions 2.5(c) and 3.8, and the fact that $b \mapsto \chi_m(b)$ is continuous.

(c)–(d) Combine Propositions 3.10–3.11.

(e) Use Lemma 3.7 and (3.58).

**4. Survival versus extinction.** In Section 4.1 we prove Theorem 1.4. In Section 4.2 we argue in favor of Conjecture 1.8.

4.1. *Proof of Theorem* 1.4. In Section 4.1.1 we introduce the Palm distribution of our process $X = (X(t))_{t\ge 0}$. In Section 4.1.2 we use the Palm distribution to identify $b_*$. In Section 4.1.3 we prove Theorem 1.4.

4.1.1. *Palm distribution and its stochastic representation.* By size-biasing our process

$$(4.1)\quad X = (X(t))_{t\ge 0}$$

at the origin at time $T$, we obtain a process

$$(4.2)\quad \widehat{X}^T = (\widehat{X}^T(t))_{0\le t\le T}$$

defined by

$$(4.3)\quad dP((\widehat{X}^T(t))_{0\le t\le T} \in \cdot) = \frac{1}{\Theta} X_0(T)\, dP((X(t))_{0\le t\le T} \in \cdot).$$

In what follows we construct a stochastic representation of $\widehat{X}^T$ in terms of a random walk and an SSDE, leading to a stochastic representation of the Palm distribution of $X$, that is, the law of $\widehat{X}^T(T)$.

*Step* 1. Fix $T > 0$. Let $\zeta = (\zeta(t))_{t\geq 0}$ be the random walk on $I$ with transition kernel $a(\cdot,\cdot)$ and jump rate 1, starting at the origin and independent of $W$. Given $\zeta$, let

$$\widehat{X}^{\zeta,T}(t) = \{\widehat{X}_i^{\zeta,T}(t)\}_{i\in I} \tag{4.4}$$

be the solution of the SSDE

$$\begin{aligned} d\widehat{X}_i^{\zeta,T}(t) &= \sum_{j\in I} a(i,j)[\widehat{X}_j^{\zeta,T}(t) - \widehat{X}_i^{\zeta,T}(t)]\,dt + \sqrt{b\widehat{X}_i^{\zeta,T}(t)^2}\,dW_i(t) \\ &\quad + b\widehat{X}_i^{\zeta,T}(t)\,1_{\{\zeta(T-t)=i\}}\,dt, \qquad i\in I, t\geq 0. \end{aligned} \tag{4.5}$$

Here, the difference with the SSDE in (1.2) is the presence of the last term, which produces a $\zeta$-dependent random potential. Like (1.2), if $\widehat{X}_i^{\zeta,T}(0)\in\mathcal{E}_1$ [the space of configurations defined in (1.4)], then (4.5) has a unique strong solution on $\mathcal{E}_1$ with continuous paths [cf. with Theorem 1.1(a)]. Let

$$\begin{aligned} P(\cdot) &= \text{ law of } \zeta \text{ on } I^{[0,\infty)}, \\ Q^{\zeta,T}(\cdot) &= \text{ law of } \widehat{X}^{\zeta,T}(T) \text{ on } [0,\infty)^I \text{ for given } \zeta. \end{aligned} \tag{4.6}$$

Then, as we prove in Proposition 4.2 below,

$$Q^T(\cdot) = \int Q^{\zeta,T}(\cdot)\,P(d\zeta) = \text{ law of } \widehat{X}^T(T) \text{ on } [0,\infty)^I \tag{4.7}$$

is the Palm distribution of $X$.

Similarly as in Lemma 2.1, we have a representation formula:

LEMMA 4.1. *Given a realization of the random walk $\zeta$, the process $(\widehat{X}^{\zeta,T}(T))_{T\geq 0}$ starting from $\widehat{X}^{\zeta,T}(0) = \underline{\Theta}$ can be represented as the following functional of the Brownian motions:*

$$\begin{aligned} \widehat{X}_i^{\zeta,T}(T) = \Theta e^{-(1/2)bT} E_i^{\xi}\Biggl(\exp\Biggl[\sqrt{b}\int_0^T \sum_{j\in I} 1_{\{\xi(T-s)=j\}}\,dW_j(s) \\ + b\int_0^T 1_{\{\xi(T-s)=\zeta(T-s)\}}\,ds\Biggr]\Biggr), \end{aligned} \tag{4.8}$$

*where $\xi = (\xi(t))_{t\geq 0}$ is the random walk on $I$ with transition kernel $a(\cdot,\cdot)$ and jump rate 1, and the expectation is over $\xi$ conditioned on $\xi(0) = i$ ($\xi$ and $W$ are independent).*

PROOF. The proof uses Itô-calculus and is similar to that of Lemma 2.1. □

*Step* 2. The link between $X(T)$ and $\widehat{X}^T(T)$ is given by the following identity, which gives an equivalent characterization of the law defined in (4.3):

PROPOSITION 4.2. *For all measurable finite $f:[0,\infty)^I \to \mathbb{R}$ that depend on finitely many coordinates,*

$$E^{\zeta,W}(f(\widehat{X}^{\zeta,T}(T))) = \frac{1}{\Theta} E^W(X_0(T)f(X(T))). \tag{4.9}$$

PROOF. Write, using Lemma 2.1,

$$\begin{aligned}
&\frac{1}{\Theta} E^W(X_0(T)f(X(T)))\\
&\quad = e^{-(1/2)bT} E^W\left(E_0^\zeta\left(\exp\left[\sqrt{b}\int_0^T \sum_{j\in I} 1_{\{\zeta(T-s)=j\}}\, dW_j(s)\right]\right) f(X(T))\right)\\
&\quad = E_0^\zeta\Biggl(E^W\Biggl(\exp\Biggl[\sqrt{b}\int_0^T \sum_{j\in I} 1_{\{\zeta(T-s)=j\}}\, dW_j(s)\\
&\qquad\qquad - \frac{1}{2} b\int_0^T \sum_{j\in I} 1_{\{\zeta(T-s)=j\}}\, ds\Biggr]\Biggr) f(X(T))\Biggr).
\end{aligned} \tag{4.10}$$

Define

$$\widehat{W}_i^\zeta(t) = W_i(t) - \sqrt{b}\int_0^t 1_{\{\zeta(T-s)=i\}}\, ds, \qquad i \in I,\ t\in[0,T]. \tag{4.11}$$

By Girsanov's formula, these are independent standard Brownian motions under the transformed path law $P^{\widehat{W}^\zeta}$ defined by

$$\begin{aligned}
\frac{dP^{\widehat{W}^\zeta}}{dP^W}(\omega_{[0,T]}) = \exp\Biggl[&\sqrt{b}\int_0^T \sum_{j\in I} 1_{\{\zeta(T-s)=j\}}\, d\omega_j(s)\\
&- \frac{1}{2} b\int_0^T \sum_{j\in I} 1_{\{\zeta(T-s)=j\}}\, ds\Biggr].
\end{aligned} \tag{4.12}$$

Now, let $\widehat{X}^{\zeta,T}$ be the unique solution of the SSDE

$$\begin{aligned}
d\widehat{X}_i^{\zeta,T}(t) = {}& \sum_{j\in I} a(i,j)[\widehat{X}_j^{\zeta,T}(t) - \widehat{X}_i^{\zeta,T}(t)]\, dt\\
&+ \sqrt{b\widehat{X}_i^{\zeta,T}(t)^2}\, d\widehat{W}_i^\zeta(t) + b\widehat{X}_i^{\zeta,T}(t)\, 1_{\{\zeta(T-t)=i\}}\, dt,\\
&\qquad\qquad i\in I,\ t\in[0,T],
\end{aligned} \tag{4.13}$$

with initial condition $\widehat{X}^{\zeta,T}(0) = \underline{\Theta}$. Comparing this with (1.2), we see that

$$\mathcal{L}(X) \text{ under } P^W = \mathcal{L}(\widehat{X}^{\zeta,T}) \text{ under } P^{\widehat{W}^\zeta}. \tag{4.14}$$

Combining this with (4.10) and (4.12), we obtain

$$\frac{1}{\Theta} E^W(X_0(T) f(X(T))) = E_0^\zeta(E^{\widehat{W}^\zeta}(f(\widehat{X}^{\zeta,T}(T)))). \tag{4.15}$$

The right-hand side of (4.15) is the same as the left-hand side of (4.9). □

*Step* 3. The reason why the Palm distribution is convenient is the following important fact, which relies on Proposition 4.2:

LEMMA 4.3. (a) w-$\lim_{T\to\infty} \mathcal{L}(\widehat{X}^T(T)) = \widehat{Q}^\infty$ *for some measure* $\widehat{Q}^\infty$ *on* $[0,\infty]^I$.

(b) *Either* $\widehat{Q}^\infty = \delta_{\underline{\infty}}$ *or* $\widehat{Q}^\infty$ *is a probability measure on* $[0,\infty)^I$.

(c) w-$\lim_{T\to\infty} \mathcal{L}(X(T)) = \delta_{\underline{0}}$ *if and only if* $\widehat{Q}^\infty = \delta_{\underline{\infty}}$.

PROOF. (a) Following Cox, Klenke and Perkins [11], via the duality relation in Lemma 2.2 we know that w-$\lim_{T\to\infty} \mathcal{L}(X(T)) = \nu_\Theta$ when $X(0) = \underline{\Theta}$ [recall (1.22)]. Since $E(X(T)) = \Theta$ for all $T \geq 0$, the Palm measure is defined for all $T \geq 0$ and converges to a limiting measure $\widehat{Q}^\infty$ on the compactified state space $[0,\infty]^I$. Indeed, because $X(T)$ converges in law, for each $n \in \mathbb{N}$ the restriction of $\widehat{X}^T(T)$ to $A_n = \{x \in [0,\infty]^I : x_0 \leq n\}$ converges in law too. The defect of $\widehat{Q}^\infty$ on $\bigcup_{n\in\mathbb{N}} A_n$ is the weight of $\underline{\infty}$.

(b) First we note that the mean of $\nu_\Theta$ is either $\Theta$ (survival) or 0 (extinction). To see why this is true, note the linearity of the SSDE: the solution for the initial condition $\Theta X(0)$ is $(\Theta X(T))_{T\geq 0}$. We know from [11], Theorem 2.3, that w-$\lim_{T\to\infty} \mathcal{L}(X(T)) = \nu_\Theta$ for every shift-invariant and shift-ergodic $\mathcal{L}(X(0))$ with mean $\Theta$. If starting from $X(0) = \underline{1}$ we get a limit law with mean $a$, then by linearity $E^{\nu_\Theta}(X_0) = a\Theta$. However, starting from $\mathcal{L}(X(0)) = \nu_\Theta$ we then get $E^{\nu_\Theta}(X_0) = a^2\Theta$ (since $\nu_\Theta$ is invariant). Thus, we must have $a^2 = a$, that is, $a = 0$ or 1.

If $\nu_\Theta$ has mean $\Theta$, then $X_0$ is uniformly integrable with respect to $\mathcal{L}(X(T))$, $T \geq 0$, and hence the Palm measures are tight. If, on the other hand, $\nu_\Theta = \delta_{\underline{0}}$, then $\widehat{Q}^\infty = \delta_{\underline{\infty}}$, because all the mass is eventually moving out of $A_n$ for every $n \in \mathbb{N}$. The representation of the Palm measure in (4.8) shows that the same holds for all other components as well.

(c) For a proof we refer to Kallenberg [36]. The intuition is the following. First note that, by translation invariance and irreducibility, the statement is equivalent to w-$\lim_{T\to\infty} \mathcal{L}(X_0(T)) = \delta_0$ if and only if w-$\lim_{T\to\infty} \mathcal{L}(\widehat{X}_0^T(T)) = \delta_\infty$. Next, if $X$ locally dies out, then (because the mean is preserved) it clusters, that is, it develops high peaks at rare sites. After size-biasing, as in Proposition 4.2, these peaks cause explosion, which is why $\widehat{X}$ diverges. Conversely, if $X$ locally survives, then it does not cluster. The size-biasing therefore does not cause explosion, which is why $\widehat{X}$ does not diverge. □

From Lemma 4.3 we obtain the following criterion for survival versus extinction of $X$. Local survival means that $\text{w-lim}_{T\to\infty}\mathcal{L}(X(T)) = \nu_\Theta \neq \delta_{\underline{0}}$.

LEMMA 4.4. *Let $P^T(\cdot)$ denote the law of $X(T)$ on $[0,\infty)^I$ [recall (4.6) and (4.7)]. Then*

$$
\begin{aligned}
(P^T)_{T\geq 0} \text{ locally survives} \quad &\Longleftrightarrow \quad (Q^T)_{T\geq 0} \text{ is tight} \\
&\Longleftrightarrow \quad (Q^{\zeta,T})_{T\geq 0} \text{ is tight } \zeta\text{-a.s.}
\end{aligned}
\tag{4.16}
$$

PROOF. The first equivalence in (4.16) uses Lemma 4.3. The second equivalence in (4.16) is trivial. □

In Section 4.1.2 we use Lemma 4.4 to identify $b_*$.

4.1.2. *Identification of $b_*$.* Abbreviate

$$
M^\xi(T) = \Theta e^{-(1/2)bT} \exp\left[\sqrt{b}\int_0^T \sum_{j\in I} 1_{\{\xi(T-s)=j\}}\, dW_j(s)\right].
\tag{4.17}
$$

This is a martingale in $W$ for every fixed $\xi$, with $M^\xi(0) = \Theta$. In terms of this quantity, the representation formulas for $X$ [in (2.2)] and $\widehat{X}$ [in (4.8)] read

$$
\begin{aligned}
X_i(T) &= E_i^\xi(M^\xi(T)),\\
\widehat{X}_i^{\zeta,T}(T) &= E_i^\xi(M^\xi(T)e^{bT^{\xi,\zeta}(T)}),
\end{aligned}
\tag{4.18}
$$

where $T^{\xi,\zeta}(T)$ is the collision local time of $\xi$ and $\zeta$ up to time $T$.

PROPOSITION 4.5. *$X$ survives for $0 < b < b_*$ and locally dies out for $b > b_*$, where*

$$
b_* = \sup\{b > 0 : (Q^{\zeta,T})_{T\geq 0} \text{ is tight } \zeta\text{-a.s.}\}.
\tag{4.19}
$$

PROOF. This follows from (4.6) and (4.16). Note that (4.8) shows that tightness of $(Q^{\zeta,T})_{T\geq 0}$ is measurable with respect to the tail sigma-algebra of $\zeta$, which is trivial. □

Equation (4.19) is to be compared with the formulas for $b_m$ and $b_{**}$ given in (1.49) and (1.51). However, (4.19) is more difficult to analyze. For one, since $b$ appears with both signs in (4.17), it is not a priori obvious that $b_*$ defines a unique transition: for this we need to show that if $(Q^{\zeta,T})_{T\geq 0}$ is tight $\zeta$-a.s. for some $b > 0$, then the same is true for all smaller values of $b$. Fortunately, the latter property can be shown to hold with the help of a coupling technique put forward in Cox, Fleischmann and Greven [9]. There it is shown that, for processes with fixed mean ($\Theta$ in our case), "more noise causes the process to be more spread out and hence to be more prone to extinction." More precisely, it is proved that two systems of the type in (1.8), with diffusion functions $g_1, g_2$ satisfying $g_1 \geq g_2$, have the property

$$E^{X,g_1}(e^{-\lambda X_0(t)}) \geq E^{X,g_2}(e^{-\lambda X_0(t)}) \qquad \forall \lambda > 0, \tag{4.20}$$

where $E^{X,g}$ is expectation over $X$ driven by $g$. Therefore, if the right-hand side converges to 1 as $t \to \infty$, then so does the left-hand side and, conversely, if the left-hand side remains bounded away from 1 as $t \to \infty$, then so does the right-hand side. Applying this to $g_1(x) = b_1 x^2$ and $g_2(x) = b_2 x^2$, we get the required monotonicity.

Our next observation is the following:

LEMMA 4.6. $b_* \geq b_{**}$.

PROOF. Take the expectation over $W$ in the second line of (4.18), use Fubini and the fact that (4.17) is a martingale in $W$, to obtain

$$E^W(\widehat{X}_i^{\zeta,T}(T)) = E_i^{\xi}(e^{bT^{\xi,\zeta}(T)}). \tag{4.21}$$

Let $b < b_{**}$. Then, according to (1.51), the right-hand side is $\zeta$-a.s. bounded in $T$. Consequently, $(\widehat{X}_i^{\zeta,T}(T))_{T\geq 0}$ is a tight family of random variables $\zeta$-a.s. Thus, by (4.6), $(Q^{\zeta,T})_{T\geq 0}$ is tight $\zeta$-a.s. Therefore, by Lemma 4.4, $b \leq b_*$. □

In Section 4.2 we will give an argument suggesting that $b_{**} > b_2$. Together with Lemma 4.6 this would imply that $b_* > b_2$, which is the claim in Conjecture 1.8.

4.1.3. *Proof of Theorem* 1.4. We look at the corresponding parts of Theorem 1.3.

(a)–(b) Combine Proposition 2.5(a) with Proposition 4.5. To get that $\nu_\Theta$ has mean $\Theta$, see the proof of Lemma 4.3.

(c)–(d) These follow from Lemma 2.2 along the lines of Cox, Klenke and Perkins [11]. The key property is the self-duality expressed in Lemma 2.2. Indeed, self-duality translates convergence of the process into convergence of the process starting from a configuration with finite mass. For the latter, martingale convergence arguments can be applied to get (d). After that, (c) follows from (b) and (d).

(e) This follows from (d) and the representation formula in Lemma 2.1, which show that $\nu_\Theta$ tends to $\delta_{\underline{\infty}}$ as $\Theta \to \infty$ because of the presence of the prefactor $\Theta$ in (2.2).

4.2. *Evidence in favor of Conjecture 1.8.* In this section we show how the conjecture follows from a certain quenched large deviation principle for a renewal process in a random environment. This large deviation principle is the topic of forthcoming work [5].

*Step* 1. According to (1.50) and (1.51), we have

$$b_{**} = \sup\left\{b > 0 : \mathrm{E}^{\xi'}\left(\exp\left[b\int_0^\infty 1_{\{\xi(t)=\xi'(t)\}}\,dt\right]\right) < \infty \ \xi\text{-a.s.}\right\}. \tag{4.22}$$

Expand the expectation in powers of $b$ to write

$$\begin{aligned}
&E^{\xi'}\left(\exp\left[b\int_0^\infty 1_{\{\xi(t)=\xi'(t)\}}\,dt\right]\right)\\
&\quad = \sum_{N=0}^\infty b^N \int_0^\infty ds_1 \int_{s_1}^\infty ds_2 \cdots \int_{s_{N-1}}^\infty ds_N\\
&\qquad\qquad \times P^{\xi'}(\xi(s_1)=\xi'(s_1),\ldots,\xi(s_N)=\xi'(s_N))\\
&\quad = \sum_{n=0}^\infty \left(\frac{b}{b_2}\right)^n K_n(\xi)
\end{aligned} \tag{4.23}$$

with

$$\begin{aligned}
K_n(\xi) = (b_2)^n &\int_0^\infty ds_1 \int_{s_1}^\infty ds_2 \cdots \int_{s_{n-1}}^\infty ds_n\\
&\times P_{s_1}(0,\xi(s_1))P_{s_2-s_1}(\xi(s_1),\xi(s_2)) \times \cdots\\
&\times P_{s_n-s_{n-1}}(\xi(s_{n-1}),\xi(s_n)),
\end{aligned} \tag{4.24}$$

where $P_t(i,j) = P^{\xi'}(\xi'(t)=j \mid \xi'(0)=i)$, $i,j \in I$, is the transition kernel of the random walk (recall Section 2.5). An easy computation shows that

$$E^{\xi}(K_n(\xi)) = 1 \qquad \forall\, n \in \mathbb{N}. \tag{4.25}$$

It therefore follows that

$$\lim_{n\to\infty} \frac{1}{n}\log K_n(\xi) = C \le 0 \qquad \xi\text{-a.s.}, \tag{4.26}$$

where the fact that the limit is $\xi$-a.s. constant is an immediate consequence of the fact that the exponential growth rate of $K_n(\xi)$ is measurable with respect to the tail sigma-algebra of $\xi$, which is trivial. By (4.22)–(4.23), to get $b_{**} > b_2$ it suffices to show that $C < 0$.

*Step* 2. Let $(\tau_k)_{k\in\mathbb{N}}$ be i.i.d. positive random variables with probability law

$$P(\tau_1 \in dt) = \frac{P_{2t}(0,0)}{(1/2)G(0,0)}\,dt. \tag{4.27}$$

Since $b_2 = 2/G(0,0)$ by (1.20), this gives

$$K_n(\xi) = E^{\tau_1,\dots,\tau_n}\left(\prod_{k=1}^{n} \frac{P_{\tau_k}(0,\Delta_k)}{P_{2\tau_k}(0,0)}\right), \tag{4.28}$$

where

$$\Delta_k = \xi\left(\sum_{m=1}^{k} \tau_m\right) - \xi\left(\sum_{m=1}^{k-1} \tau_m\right), \qquad k \in \mathbb{N}, \tag{4.29}$$

with $\tau_0 = 0$. The expectation in (4.28) seems well suited for a large deviation analysis, but the problem is that the $\Delta_k$ are correlated because $\xi$ is *fixed*.

*Step* 3. Abbreviate

$$T_k = \sum_{l=1}^{k} \tau_l, \qquad Y_k = \xi_{[T_{k-1},T_k)}, \qquad k \in \mathbb{N}. \tag{4.30}$$

Then $(Y_k)_{k\in\mathbb{N}}$ is a random process taking values in $\Omega = \bigcup_{t\geq 0} D([0,t),I)$, with $D([0,t),I)$ the Skorohod space of paths of length $t$. For $z = z_{[0,t)} \in \Omega$ we abbreviate

$$\Delta z = z(t) - z(0), \qquad |z| = t - 0 = t. \tag{4.31}$$

Introduce the empirical process

$$R_n = \frac{1}{n}\sum_{k=1}^{n} \delta_{(Y_k,\dots,Y_{k+n-1})^{per}}, \qquad n \in \mathbb{N}, \tag{4.32}$$

where the upper index *per* means that the $n$-sequence is periodically repeated. For each $k \in \mathbb{N}$, we have $Y_k \in \Omega$. Consequently, $R_n$ is a random element of $\mathcal{P}(\Omega^{\mathbb{N}})$, the set of probability measures on $\Omega^{\mathbb{N}}$, with the randomness coming from $(\tau_k)_{k\in\mathbb{N}}$ ($\xi$ being fixed). In terms of $R_n$ we may write

$$K_n(\xi) = E^{\tau_1,\dots,\tau_n}\left(\exp\left[n\int_{\Omega} (\pi_1 R_n)(dy) \log \frac{P_{|y|}(0,\Sigma y)}{P_{2|y|}(0,0)}\right]\right), \tag{4.33}$$

where $\pi_1 R_n$ is the projection of $R_n$ onto the first coordinate.

*Step* 4. In [5] it is shown that $(R_n)_{n\in\mathbb{N}}$ satisfies a quenched (with respect to $\xi$) large deviation principle on $\mathcal{P}(\Omega^{\mathbb{N}})$ with some rate function $Q \mapsto I(Q)$. This rate function turns out to be a sum of two terms, the first being the specific relative entropy of $Q$ with respect to $Q^0$, where

$$Q^0 = \mathcal{L}((Y_k)_{k\in\mathbb{N}}), \tag{4.34}$$

the second being a specific relative entropy associated with "concatenation and randomization" of $(Y_k)_{k\in\mathbb{N}}$. By applying Varadhan's lemma to (4.33), we therefore obtain

$$\lim_{n\to\infty} \frac{1}{n} \log K_n(\xi) = C, \qquad \xi\text{-a.s.} \tag{4.35}$$

with

$$C = \sup_{Q\in\mathcal{P}^{\mathrm{stat}}(\Omega^{\mathbb{N}})} \left[ \int_\Omega (\pi_1 Q)(dy) \log \frac{P_{|y|}(0,\Sigma y)}{P_{2|y|}(0,0)} - I(Q) \right]. \tag{4.36}$$

We note that Varadhan's lemma applies because $(t,i) \mapsto \log P_t(0,i)/P_{2t}(0,0)$ is bounded from above and tending to $-\infty$ as $\|i\| \to \infty$ (see Deuschel and Stroock [21], Theorem 2.1.10).

*Step* 5. Note that $Q^0$ is a product measure:

$$Q^0 = (q^0)^{\mathbb{N}} \qquad \text{with } q^0(dz) = q^0(d\tau, d\xi_{[0,\tau)}) = \nu(d\tau)\mu(d\xi_{[0,\tau)}). \tag{4.37}$$

Since the integral in (4.36) depends on $\pi_1 Q$ only, by putting

$$\inf_{\substack{Q\in\mathcal{P}^{\mathrm{stat}}(\Omega^{\mathbb{N}})\\ \pi_1 Q = q}} I(Q) = i(q), \tag{4.38}$$

we can reduce the variational expression in (4.36) to

$$C = \sup_{q\in\mathcal{P}(\Omega)} \left[ \int_\Omega q(dy) \log \frac{P_{|y|}(0,\Sigma y)}{P_{2|y|}(0,0)} - i(q) \right]. \tag{4.39}$$

This formula is the key to proving that $C < 0$. Indeed, $i(q)$ turns out to be a sum of two terms, one being $h(q|q^0)$ (the relative entropy of $q$ with respect to $q^0$), the other being a relative entropy term associated with concatenation and randomization. Without the second term, the supremum in (4.39) would be attained at

$$q(dz) = q(d\tau, d\xi_{[0,\tau)}) = b_2 P_\tau(0,\xi(\tau))\, d\tau\, \mu(d\xi_{[0,\tau)}) \tag{4.40}$$

and would be equal to 0. Indeed, this corresponds to the annealed upper bound in (4.25) and (4.26). With the second term, however, it is $< 0$. See [5] for further details.

**Acknowledgments.** The authors thank Matthias Birkner (Berlin), Ted Cox (Syracuse) and Frank Redig (Leiden) for valuable input during the preparation of this work, an anonymous referee for critical comments, and the staff of EURANDOM for multiple hospitality.

Mathematisches Institut
Universität Erlangen–Nürnberg
Bismarckstrasse $1\frac{1}{2}$
D-91504 Erlangen
Germany
E-mail: greven@mi.uni-erlangen.de

Mathematical Institute
Leiden University
P.O. Box 9512
2300 RA Leiden
and
EURANDOM
P.O. Box 513
5600 MB Leiden
The Netherlands
E-mail: denholla@math.leidenuniv.nl